%% file: promotion.tex
\documentclass{amsart}
\usepackage{verbatim,amsmath,amssymb}
\usepackage{ifthen}
\usepackage[all]{xy}
\usepackage[usenames,svgnames]{xcolor}
\usepackage{tikz}
\usetikzlibrary{matrix}
\usepackage{tabmac}
\newcommand{\tab}{\tableau[sbY]}

\newtheorem{theorem}{Theorem}
\newtheorem{proposition}[theorem]{Proposition}
\newtheorem{lemma}[theorem]{Lemma}

\newtheorem{conjecture}[theorem]{Conjecture}

\theoremstyle{definition}
\newtheorem{definition}[theorem]{Definition}
\newtheorem{remark}[theorem]{Remark}
\newtheorem{remarks}[theorem]{Remarks}
\newtheorem{example}[theorem]{Example}

\numberwithin{equation}{section}
\numberwithin{theorem}{section}

\newboolean{draft}
\setboolean{draft}{true}
\ifdraft
\newcommand{\TODO}[2][To do: ]{\textcolor{red}{\textbf{#1#2}}}
\else
\newcommand{\TODO}[2][]{}
\fi

\newcommand{\id}{\mathrm{id}}
\newcommand{\Pcl}{P_{\mathrm{cl}}}
\newcommand{\pr}{\mathrm{pr}}
\newcommand{\prA}{\mathfrak{pr}}
\newcommand{\prto}{\overset{\pr}{\longrightarrow}}
\newcommand{\wt}{\mathrm{wt}}
\newcommand{\ZZ}{\mathbb{Z}}
\newcommand{\fail}{\emptyset}

\usepackage{verbatim}
\makeatletter
\def\Mexin@processline{>> \the\verbatim@line\par}
\newenvironment{Mexin} {\vspace{-.5ex}\verbatim\small\addtolength\parskip{-.5ex}\let\verbatim@processline=\Mexin@processline}{\endverbatim}
\newenvironment{Mexout}{\vspace{-.5ex}\verbatim\small\addtolength\parskip{-.9ex}}{\endverbatim}
\makeatother

\makeatletter

\newskip\@bigflushglue \@bigflushglue = -100pt plus 1fil

\def\bigcentering{\let\\\@centercr\rightskip\@bigflushglue%
\leftskip\@bigflushglue
\parindent\z@\parfillskip\z@skip}

\makeatother

\begin{document}

\title[Promotion operator on tensor products of crystals]{On the uniqueness of promotion operators on 
tensor products of type $A$ crystals}

\author[J. Bandlow]{Jason Bandlow}
\address{Department of Mathematics, University of Pennsylvania, 
209 South 33rd Street, Philadelphia, PA 19104-6395, U.S.A.} 
\email{jbandlow@math.upenn.edu}
\urladdr{http://www.math.upenn.edu/\~{}jbandlow}

\author[A. Schilling]{Anne Schilling}
\address{Department of Mathematics, University of California, One Shields
Avenue, Davis, CA 95616-8633, U.S.A.}
\email{anne@math.ucdavis.edu}
\urladdr{http://www.math.ucdavis.edu/\~{}anne}

\author[N. Thi\'ery]{Nicolas M. Thi\'ery}
\address{Univ Paris-Sud, Laboratoire de Math\'ematiques d'Orsay,
  Orsay, F-91405; CNRS, Orsay, F-91405, France}
\address{\textit{Current address}: Department of Mathematics, University of California, One
 Shields Avenue, Davis, CA 95616, U.S.A.}
\email{Nicolas.Thiery@u-psud.fr}
\urladdr{http://Nicolas.Thiery.name}

\thanks{\textit{Date:} June 2008; updated May 2009}


\begin{abstract}
    The affine Dynkin diagram of type $A_n^{(1)}$ has a cyclic symmetry. The analogue of this Dynkin
    diagram automorphism on the level of crystals is called a promotion operator. In this paper we 
    show that the only irreducible type $A_n$ crystals which admit a 
    promotion operator are the highest weight crystals indexed by rectangles. In addition we prove
    that on the tensor product of two type $A_n$ crystals labeled by
    rectangles, there is a single connected
    promotion operator. We conjecture this to be true for an arbitrary number of tensor factors. Our 
    results are in agreement with Kashiwara's conjecture that all `good' affine crystals are tensor 
    products of Kirillov-Reshetikhin crystals.
\end{abstract}  

\maketitle

\section{Introduction}
\label{section.introduction}
The Dynkin diagram of affine type $A_n^{(1)}$ has a cyclic symmetry generated by
the map $i\mapsto i+1 \pmod{n+1}$. The promotion operator is the analogue of this
Dynkin diagram automorphism on the level of crystals.
Crystals were introduced by Kashiwara~\cite{Kashiwara.1994} to give a combinatorial 
description of the structure of modules over the universal enveloping algebra
$U_q(\mathfrak g)$ when $q$ tends to zero.
In short, a crystal is a non-empty set $B$ endowed with raising and
lowering crystal operators $e_i$ and $f_i$ indexed by the nodes of the
Dynkin diagram $i\in I$, as well as a weight function $\wt$. It can be
depicted as an edge-colored directed graph with elements of $B$ as
vertices and $i$-arrows given by $f_i$. In type $A_n$, the highest weight
crystal $B(\lambda)$ of highest weight $\lambda$ is the set of all
semi-standard Young tableaux of shape $\lambda$ (see for
example~\cite{Kashiwara.Nakashima.1994,LLT.1995}) with weight function
given by the content of tableaux.

\begin{definition} \label{definition.promotion}
A \emph{promotion operator} $\pr$ on a crystal $B$ of type $A_n$ is an operator $\pr:B\to B$
such that:
  \begin{enumerate}
    \item \label{pt.content} $\pr$ shifts the content: If $\wt(b)=(w_1,\ldots,w_{n+1})$ is the content 
      of the crystal element $b\in B$, then $\wt(\pr(b))=(w_{n+1},w_1,\ldots,w_n)$;
    \item \label{pt.order} Promotion has order $n+1$: $\pr^{n+1} = \id$;
    \item \label{pt.e} $\pr \circ e_i = e_{i+1} \circ  \pr$ and $\pr \circ f_i = f_{i+1} \circ  \pr$ for 
      $i\in \{1,2,\ldots,n-1\}$.
  \end{enumerate}
If condition~\eqref{pt.order} is not satisfied, but $\pr$ is still
bijective, then $\pr$ is a \emph{weak} promotion operator.
\end{definition}

Given a (weak) promotion operator on a crystal $B$ of type $A_n$, one can define an 
associated (weak) \emph{affine crystal} by setting
\begin{equation}
  \label{eq.affineCrystal}
	e_0 := \pr^{-1} \circ e_1 \circ \pr \quad \text{and} \quad f_0 := \pr^{-1} \circ f_1 \circ \pr.
\end{equation}
A promotion operator $\pr$ is called \emph{connected} if the resulting affine crystal $B$ is 
connected (as a graph).  Two promotion operators are called \emph{isomorphic} if the resulting
affine crystals are isomorphic.

Our aim is the classification of all affine crystals that are
associated to a promotion operator on a tensor product of 
highest weight crystals $B(\lambda)$ of type $A_n$.

Sch\"utzenberger~\cite{Schuetzenberger.1972} introduced a weak
promotion operator $\prA$ on tableaux using jeu-de-taquin (see
Section~\ref{subsection.promotion}). It turns out that $\prA$ is the
unique weak promotion operator on $B(\lambda)$; furthermore, $\prA$ is
a promotion operator if and only if $\lambda$ is a rectangle
(cf. Proposition~\ref{proposition.prA} which is based on results by
Haiman~\cite{Haiman.1992} and Shimozono~\cite{Shimozono.2002}).

Let us denote by $\omega_1,\dots,\omega_n$ the fundamental weights of
type $A_n$. One can identify the rectangle partition $\lambda:=(s^r)$
of height $r$ and width $s$ with the weight $s\omega_r$. We henceforth
call $\prA$ on $B(s\omega_r)$ the \emph{canonical promotion
  operator}. It can be extended to tensor products $B(s_1\omega_{r_1})
\otimes \cdots \otimes B(s_\ell \omega_{r_\ell})$ indexed by
rectangles by setting $\prA(b_1\otimes \cdots \otimes b_\ell) :=
\prA(b_1) \otimes \cdots \otimes \prA(b_\ell)$.  Let $B$ be a crystal
with an isomorphism $\Psi$ to a direct sum of tensor products of
highest weight crystals indexed by rectangles. A promotion
operator is \emph{induced by $\Psi$} if it is of the form
$\Psi^{-1}\circ\prA\circ \Psi$, where $\prA$ is the canonical promotion on each
summand. Note that throughout the paper, all tensor factors are written in reverse direction
compared to Kashiwara's conventions, which is more compatible with operations on tableaux.

The main result of this paper is the following theorem.
\begin{theorem}
  \label{theorem.main}
  Let $B=B(s'\omega_{r'}) \otimes B(s\omega_r)$ be the tensor product
  of two classical highest weight crystals of type $A_n$ with $n\ge
  2$, labeled by rectangles. If $(s,r)\neq (s',r')$, there is a unique
  promotion operator $\pr=\prA$. If $(s,r)=(s',r')$, there are two
  promotion operators: The canonical one $\pr=\prA$ which is connected
  and the one induced by $\Psi$ (with $\Psi$ as defined
  in~\eqref{equation.isomorphism}) which is disconnected.
\end{theorem}

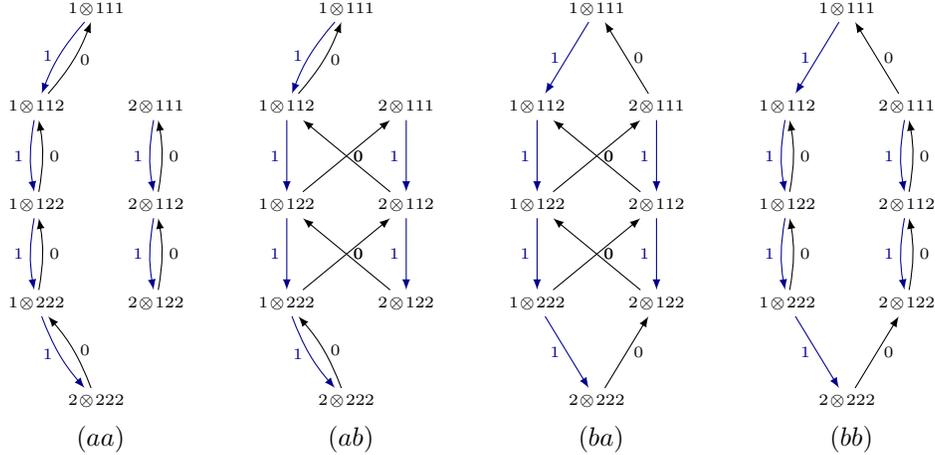
\begin{figure}
  \begin{bigcenter}
    \begin{tabular}{c@{\qquad}c@{\qquad}c@{\qquad}c}
      \input{figure-promotions-B11B13-1} &
      \input{figure-promotions-B11B13-2} &
      \input{figure-promotions-B11B13-3} &
      \input{figure-promotions-B11B13-4} \\
      $(aa)$ & $(ab)$ & $(ba)$ & $(bb)$
    \end{tabular}
  \end{bigcenter}
  \caption{The four affine crystals associated to the classical crystal
    $B(\omega_1)\otimes B(3\omega_1)$ for type $A_1$. The affine crystal
    $B^{1,1}\otimes B^{3,1}$ corresponds to $(bb)$. The others are not
    'good' crystals (see Definition~\ref{definition.good}):  (aa) is not connected,
    (ab) is not simple, and (ba) does not satisfy
    the convexity condition on string lengths.
  }
  \label{fig:promotionsForA11}
\end{figure}

\begin{remark} \label{remark.counterexample}
  As illustrated in Figure~\ref{fig:promotionsForA11},
  Theorem~\ref{theorem.main} does not hold for $n=1$.
  Only (bb) yields a 'good' crystal according to the combinatorial
  Definition~\ref{definition.good}. It would be interesting to determine whether
  (ab) and (ba) correspond to crystals for $U_q'(\widehat{\mathfrak{sl}}_2)$-modules.
\end{remark}

As suggested by further evidence discussed in
Section~\ref{section.evidence}, we expect this result to carry over to
any number of tensor factors.
\begin{conjecture}
  \label{conjecture.promotion} 
  Let $B:=B(\lambda^{1}) \otimes \cdots \otimes B(\lambda^{\ell})$ be
  a tensor product of classical highest weight crystals of type $A_n$
  with $n\ge 2$.  Then, any promotion operator is induced by an
  isomorphism $\Psi$ from $B$ to some direct sum of tensor products of
  classical highest weight crystals of rectangular shape.

  Furthermore, there exists a connected promotion operator if and only
  if $\lambda^{1},\dots,\lambda^{\ell}$ are rectangles, and this
  operator is $\prA$ up to isomorphism.
\end{conjecture}

As shown by Shimozono~\cite{Shimozono.2002}, the affine crystal constructed from
$B(s\omega_r)$ using the promotion operator $\prA$ is isomorphic to the Kirillov-Reshetikhin
crystal $B^{r,s}$ of type $A_n^{(1)}$. Kirillov-Reshetikhin crystals $B^{r,s}$ form a special class of
finite dimensional affine crystals, indexed by a node $r$ of the classical Dynkin diagram and 
a positive integer $s$.
Finite-dimensional affine $U_q'(\mathfrak{g})$-crystals have been used extensively in
the study of exactly solvable lattice models in statistical mechanics. 
It has recently been proven~\cite{KKMMNN.1992,Okado_Schilling.2008} that 
(for nonexceptional types) the Kirillov-Reshetikhin module $W(s\omega_r)$, labeled by a positive 
multiple of the fundamental weight $\omega_r$, has a crystal basis called the Kirillov-Reshetikhin
crystal $B^{r,s}$. Kashiwara conjectured (see Conjecture~\ref{conjecture.Kashiwara}) that any
`good' affine finite crystal is the tensor product of Kirillov-Reshetikhin crystals.

Note that Theorem~\ref{theorem.main} and Conjecture~\ref{conjecture.promotion} are in
agreement with Kashiwara's Conjecture~\ref{conjecture.Kashiwara}.
Namely, if one can assume that every `good' affine crystal for type $A_n^{(1)}$ comes 
from a promotion operator, then Theorem~\ref{theorem.main} and Conjecture~\ref{conjecture.promotion}
imply that any crystal with underlying classical crystal being a tensor product is a tensor product of 
Kirillov-Reshetikhin crystals.

Promotion operators have appeared in other contexts as well. Promotion has been studied
by Rhoades et al.~\cite{Rhoades.2008, PPR.2008} in relation with Kazhdan-Lusztig theory 
and the cyclic sieving phenomenon. 
Hernandez~\cite{Hernandez.2008} proved $q$-character formulas for cyclic Dynkin diagrams in 
the context of toroidal algebras. He studies a ring morphism $R$ which is related to the promotion 
of the Dynkin diagram. Since $q$-characters are expected to be related to crystal theory, this is 
another occurrence of the promotion operator.
Theorem~\ref{theorem.main} is also a first step in defining an affine crystal on rigged 
configurations. There exists a bijection between tuples of rectangular tableaux and rigged 
configurations~\cite{Kerov.Kirillov.Reshetikhin.1986,Kirillov.Schilling.Shimozono.2002,
Deka.Schilling.2006}. A classical crystal on rigged configurations was defined
and a weak promotion operator was conjectured in~\cite{Schilling.2006}. It remains to
prove that this weak promotion operator has the correct order.

This paper is organized as follows. In
Section~\ref{section.typeACrystals} crystal theory for type $A_n$ is
reviewed, some basic properties of promotion operators are stated
which are used later, and Kashiwara's conjecture is stated. In
Section~\ref{section.promotion}, the Sch\"utzenberger map $\prA$ is
defined on $B(\lambda)$ using jeu-de-taquin. It is shown that it is
the only possible weak promotion operator on $B(\lambda)$, and that it
is a promotion operator on $B(\lambda)$ if and only if $\lambda$ is of
rectangular shape. Section~\ref{section.proof} is devoted to the proof
of Theorem~\ref{theorem.main} and in Section~\ref{section.evidence} we
provide evidence for Conjecture~\ref{conjecture.promotion}; in
particular, we discuss unique factorization into a product of Schur
polynomials indexed by rectangles.

\section*{Acknowledgements}
\label{section.acknowledgements}
We would like to thank Ghislain Fourier, Mark Haiman, David Hernandez, Masaki 
Kashiwara, Masato Okado, Brendon Rhoades, Mark Shimozono, and John Stembridge
for fruitful discussions. Special thanks to Stephanie van Willigenburg for pointing us to the 
references~\cite{Purbhoo.Willigenburg.2007} and~\cite{Rajan.2004}.
This research was partially supported by 
NSF grants DMS-0501101, DMS-0652641, and DMS-0652652.  We would also like
to thank the inspiring program on ``Combinatorial Representation
Theory'' held at MSRI from January through May 2008 during which time
much of this research was carried out.

The research was partially driven by computer exploration using the
open-source algebraic combinatorics package
\texttt{MuPAD-Combinat}~\cite{MuPAD-Combinat}, together with
\texttt{Sage}~\cite{Sage}. The pictures in this paper have been produced
(semi)-automatically, using \texttt{MuPAD-Combinat},
\texttt{graphviz}, \texttt{dot2tex}, and \texttt{pgf/tikz}.

\section{Review of type $A$ crystals}
\label{section.typeACrystals}

In this section, we recall some definitions and properties of type $A$ crystals, state some 
lemmas which will be used extensively in the proof of Theorem~\ref{theorem.main}, and
state Kashiwara's conjecture.  

\subsection{Type $A$ crystal operations}
Crystal graphs of integrable $U_q(\mathfrak{sl}_{n+1})$-modules can be
defined by operations on tableaux (see for
example~\cite{Kashiwara.Nakashima.1994,LLT.1995}).  Consider the type
$A_n$ Dynkin diagram with nodes indexed by $I:=\{1,\dots,n\}$. There
is a natural correspondence between dominant weights in the weight
lattice $P:=\bigoplus_{i\in I} \ZZ \omega_i$, where $\omega_i$ is the
$i$-th fundamental weight, and partitions $\lambda=(\lambda_1\ge
\lambda_2 \ge \cdots \ge \lambda_n \ge 0)$ with at most $n$
parts. Suppose $\lambda=\omega_{r_1}+\cdots+\omega_{r_k}$ is a dominant
weight. Then we can associate to $\lambda$ the partition with columns
of height $r_1,\ldots,r_k$. In particular, the fundamental weight
$s\omega_r$ is associated to the partition of rectangular shape of
width $s$ and height $r$.

The highest weight crystal $B(\lambda)$ of type $A_n$ is given by the set of all 
semi-standard tableaux of shape $\lambda$ over the alphabet $\{1,2,\ldots,n+1\}$ endowed
with maps
\begin{equation*}
\begin{split}
	e_i, f_i: &B(\lambda) \to B(\lambda)\cup \{\fail\} \quad \text{for $i\in I=\{1,2,\ldots,n\}$},\\
	\wt : &B(\lambda) \to P.
\end{split}
\end{equation*}
Throughout this paper, we use French notation for tableaux (that is,
they are weakly increasing along rows from left to right and strictly
increasing along columns from bottom to top).  The \emph{weight} of a
tableau $t$ is its content
$$\wt(t):=(m_1(t),m_2(t),\ldots,m_{n+1}(t))\,,$$ where $m_i(t)$ is the
number of letters $i$ appearing in $t$. The \emph{lowering and raising
  operators} $f_i$ and $e_i$ can be defined as follows. Consider the
\emph{row reading word} $w(t)$ of $t$; it is obtained by reading the
entries of $t$ from left to right, top to bottom. Consider the subword
of $w(t)$ consisting only of the letters $i$ and $i+1$ and associate
an open parenthesis '$)$' with each letter $i$ and a closed parenthesis
'$($' with each letter $i+1$. Successively match all parentheses. Then
$f_i$ transforms the letter $i$ that corresponds to the rightmost
unmatched parenthesis '$)$' into an $i+1$. If no such parenthesis '$)$'
exists, $f_i(t)=\fail$. Similarly, $e_i$ transforms the letter $i+1$ that
corresponds to the leftmost unmatched parenthesis '$($' into an $i$. If
no such parenthesis exists, $e_i(t)=\fail$.

For a tableau $t$, define $\varphi_i(t)=\max \{k \mid f_i^k(t)=\fail\}$ (resp. $\varepsilon_i(t)=
\max \{k \mid e_i^k(t)=\fail\}$) to be the maximal number of times $f_i$ (resp. $e_i$) can be applied to 
$t$. The quantity $\varphi_i(t)+\varepsilon_i(t)$ is the length of the \emph{$i$-string} of $t$.
Similarly, let $b_i(t)$ be the number of paired '$()$' parentheses in the algorithm for computing $f_i$ 
and $e_i$. We call this the number of \emph{$i$-brackets} in $t$.

\begin{example}
Let
\begin{equation*}
         \def\t{\tableau[sbY]}
         \let\c=\vcenter
	t = \c{\t{2,3,3|1,2,2,3}} \; .
\end{equation*}
Then $w(t)=2331223$ and
\begin{equation*}
	\def\t{\tableau[sbY]}
	\let\c=\vcenter
	f_2(t) = \c{\t{3,3,3|1,2,2,3}} \quad \text{and} \quad 
	e_2(t) = \c{\t{2,3,3|1,2,2,2}} \; .
\end{equation*}
\end{example}

\begin{definition}
  For $J \subset I=\{1,2,\ldots,n\}$, the element $b\in B$ is
  \emph{$J$-highest weight} if $e_i(b)=\fail$ for all $i\in J$. It is
  \emph{highest weight} if it is $I$-highest weight. Similarly, $b\in B$ is
  \emph{$J$-lowest weight} if $f_i(b)=\fail$ for all $\in J$.
\end{definition}

\subsection{Crystal isomorphisms}
Let $B$ and $B'$ be two crystals over the same Dynkin diagram. Then a
bijective map $\Phi:B\to B'$ is a \emph{crystal isomorphism} if for all $b\in
B$ and $i\in I$,
\begin{equation*}
	f_i \Phi(b) = \Phi(f_i b) \quad \text{and} \quad e_i \Phi(b) = \Phi(e_i b)\,,
\end{equation*} 
where by convention $\Phi(\fail)=\fail$.  More generally, let $B$ and $B'$ be
crystals over two isomorphic Dynkin diagrams $D$ and $D'$ with nodes
respectively indexed by $I$ and $I'$, and let $\tau:I\to I'$ be an
isomorphism from $D$ to $D'$.  Then $\Phi$ is a
$\tau$-\emph{twisted-isomorphism} if for all $b\in B$ and $i\in I$,
\begin{equation*}
  f_{\tau(i)} \Phi(b) = \Phi(f_i b) \quad \text{and} \quad e_{\tau(i)} \Phi(b) = \Phi(e_i b)\,.
\end{equation*}

It was proven by Stembridge~\cite{Stembridge.2001} that in the expansion of the product of 
two Schur functions indexed by rectangles, each summand $s_\lambda$ occurs with multiplicity 
zero or one. This implies in particular that in the decomposition of type $A_n$ crystals
\begin{equation} \label{eq:classicalDecomposition}
	B(s'\omega_{r'}) \otimes B(s\omega_r) \cong \bigoplus_{\lambda} B(\lambda)
\end{equation}
each irreducible component $B(\lambda)$ occurs with multiplicity at most one. Hence there
is a \emph{unique} crystal isomorphism
\begin{equation} \label{eq:unique_iso}
	B(s'\omega_{r'}) \otimes B(s\omega_r) \cong B(s\omega_r) \otimes B(s'\omega_{r'}).
\end{equation}
Recall that all tensor factors are written in reverse direction compared to Kashiwara's conventions.

For two equal rectangular tensor factors, there is a unique additional crystal isomorphism
\begin{equation} \label{equation.isomorphism}
	\Psi: B(s\omega_r)^{\otimes 2} \cong B((s-1)\omega_r) \otimes B((s+1)\omega_r) 
	\oplus B(s\omega_{r-1}) \otimes B(s\omega_{r+1}).
\end{equation}
Its existence follows from the well-known Schur function equality
\cite{Kirillov.1983, Kleber.2001}:
\begin{displaymath}
  s_{(s^r)}^2 = s_{((s-1)^{r})} s_{((s+1)^r)} + s_{(s^{r-1})} s_{(s^{r+1})}. 
\end{displaymath}
This isomorphism can be described explicitly as follows. For 
$b'\otimes b\in B(s\omega_r)^{\otimes 2}$ consider the tableau $b'.b$ given by the Schensted 
row insertion of $b$ into $b'$. By~\cite{Stembridge.2001}, there is a unique pair of tableaux 
$\tilde{b}'\otimes \tilde{b}$ either in $B((s-1)\omega_r) \otimes B((s+1)\omega_r)$ or in
$B(s\omega_{r-1}) \otimes B(s\omega_{r+1})$ such that $\tilde{b}'.\tilde{b} = b'.b$. 
Define $\Psi(b'\otimes b) = \tilde{b}' \otimes \tilde{b}$.

\begin{example}
Let
\begin{equation*}
        \def\t{\tableau[sbY]}
         \let\c=\vcenter
	b' \otimes b = \c{\t{2,3|1,2}} \otimes \c{\t{2,2|1,1}} \quad \text{so that} \quad
	b'.b = \c{\t{3|2,2,2|1,1,1,2}} \; .
\end{equation*}
Then
\begin{equation*}
        \def\t{\tableau[sbY]}
         \let\c=\vcenter
	\Psi(b' \otimes b) = \tilde{b}' \otimes \tilde{b} = \c{\t{3|2}} \otimes \c{\t{2,2,2|1,1,1}} 
	\quad \text{since} \quad
	\tilde{b}'.\tilde{b} = b'.b.
\end{equation*}
If on the other hand 
\begin{equation*}
        \def\t{\tableau[sbY]}
         \let\c=\vcenter
	b' \otimes b = \c{\t{3,3|1,2}} \otimes \c{\t{2,2|1,1}} \quad \text{then} \quad
	b'.b = \c{\t{3,3|2,2|1,1,1,2}} \; .
\end{equation*}
Hence 
\begin{equation*}
        \def\t{\tableau[sbY]}
         \let\c=\vcenter
	\Psi(b' \otimes b) = \tilde{b}' \otimes \tilde{b} = \c{\t{1,2}} \otimes \c{\t{3,3|2,2|1,1}}\; . 
\end{equation*}
\end{example}

\subsection{Duality}
For each $A_n$ crystal $B(\lambda)$ of highest weight $\lambda$, there exists a dual
crystal $B(\lambda^\complement)$, where $\lambda^\complement$ is the \emph{complement partition}
of $\lambda$ in a rectangle of height $n+1$ and width $\lambda_1$. The crystal $B(\lambda)$ and 
its dual $B(\lambda^\complement)$ are twisted-isomorphic, with $\tau(i)=n+1-i$.

\begin{proposition} \label{proposition.duality}
   The $A_n$ crystal $B=B(s_1\omega_{r_1}) \otimes \cdots \otimes B(s_\ell \omega_{r_\ell})$ 
   is twisted-isomorphic to the $A_n$ crystal $B(s_1\omega_{n+1-r_1}) \otimes \cdots \otimes 
   B(s_\ell \omega_{n+1-r_\ell})$.
\end{proposition}
\begin{proof}
 This follows from the fact that the tensor product of twisted-isomorphic crystals must be 
 twisted-isomorphic.
\end{proof}
  
 \begin{lemma}[Duality Lemma] \label{lemma.duality}
   All promotion operators on $B=B(s'\omega_{r'}) \otimes B(s\omega_r)$ of type $A_n$ are in
   one-to-one connectedness-preserving correspondence with the promotion operators on 
   $B(s'\omega_{n+1-r'}) \otimes B(s\omega_{n+1-r})$. As a consequence, to classify all promotion 
   operators on $B$, it suffices to classify them for $n \le r + r' - 1$.
 \end{lemma}
 \begin{proof}
   By Proposition~\ref{proposition.duality}, $B(s'\omega_{r'}) \otimes B(s\omega_r)$ is 
   twisted-isomorphic to $B(s'\omega_{n + 1 - r'}) \otimes B(s \omega_{n + 1 - r})$.  
   Notice that under this twisted-isomorphism $\Phi$, a promotion $\pr$ on $B$ becomes
   $\Phi\circ \pr$ and satisfies the conditions of Definition~\ref{definition.promotion} of the inverse 
   of a promotion. Hence each $\pr$ induces a promotion on the dual of $B$. It is clear
   that connectedness is preserved.
   
   Now suppose $n>r + r' -1$. Summing the heights of the dual tensor product and subtracting one, 
   we obtain
   \begin{equation*}
     (n + 1 - r') + (n + 1 - r) - 1 = 2n - (r + r') +1 > n,
   \end{equation*}
   which satisfies the condition of the lemma. Hence it suffices to classify promotion operators
   for $n\le r+r'-1$.
\end{proof}

\subsection{Properties of promotion operators}

In this section we discuss some further properties of promotion
operators. We begin with two remarks about consequences of the axioms
for a promotion operator as defined in
Definition~\ref{definition.promotion} which will be used later.  In
particular, in Remark~\ref{remark.shift2} a reformulation of the three
conditions in Definition~\ref{definition.promotion} is provided which
in practice might be easier to verify. Then we prove two Lemmas: the
Highest Weight Lemma~\ref{lemma.highestWeight} and the Two Path
Lemma~\ref{lemma.twoPaths}.

\begin{remark}
  \label{remark.shift}
  Let $\pr$ be a promotion operator. Then, $\pr^k \circ e_i = e_{i+k}
  \circ \pr^k$ whenever $i,i + k \neq 0 \pmod{n+1}$, and similarly for
  $f_i$.
\end{remark}
\begin{proof}
  Iterate condition~\eqref{pt.e} of Definition~\ref{definition.promotion}, using condition~\eqref{pt.order} 
  to go around $i=0$.
\end{proof}

\begin{remark} 
  \label{remark.shift2} 
  Let $B:=B_1\otimes\cdots\otimes B_\ell$ be a tensor product of type
  $A_n$ highest weight crystals (or more generally a crystal of type $A_n$
  with some weight space of dimension $1$; this includes the simple
  crystals of Definition~\ref{definition.simple}), and $\pr$ a weak promotion
  operator on $B$ which satisfies:
  \begin{enumerate}
  \item[(2')]\label{pt.e2}  $\pr^2 \circ e_n = e_1 \circ \pr^2$, and $\pr^2 \circ f_n = f_1 \circ \pr^2$.
  \end{enumerate}
  Assume that the associated weak affine crystal graph is connected. Then, $\pr$ is a promotion 
  operator.
\end{remark}
\begin{proof}
  We need to prove condition~\eqref{pt.order}: $\pr^{n+1}=\id$. First
  note that condition (2') together with the definition of $e_0$
  in~\eqref{eq.affineCrystal} implies that condition~\eqref{pt.e}:
  $\pr \circ e_i = e_{i+1} \circ \pr$ (with $i+1$ taken $\pmod{n+1}$)
  holds even for $i=n$. By repeated application, one obtains
  $\pr^{n+1} \circ e_i = e_i \circ \pr^{n+1}$ for all $i$ (and
  similarly for $f_i$). In other words, $\pr^{n+1}$ is an automorphism
  of the weak affine crystal graph.

  We now check that such an automorphism has to be trivial. First note
  that it preserves classical weights. For all $1\le j\le \ell$, let $u_j$ be the
  highest vector of $B_j$. Then, $u:=u_1\otimes\cdots\otimes u_\ell$
  is the unique element of $B$ of weight $\wt(u_1)+\dots+\wt(u_\ell)$,
  and therefore is fixed by $\pr^{n+1}$.  Take finally any $v\in
  B$. By the connectivity assumption $v=F(u)$, where $F$ is
  some concatenation of crystal operators.  Therefore, $\pr^{n+1}(v) =
  \pr^{n+1} \circ F(u) = F (\pr^{n+1}(u)) = F(u)=v$.
\end{proof}

For the remainder of this section $B$ is a crystal of type $A_n$ on which a promotion operator
$\pr$ is defined. Recall that for $J \subset \{1,2,\ldots,n\}$, the element $b\in B$ is 
\emph{$J$-highest weight} if $e_i(b)=\fail$ for all $i\in J$.
  
\begin{lemma}[Highest Weight Lemma] \label{lemma.highestWeight}
  If $\pr(b)$ is known for all $\{1,2,\ldots,n-1\}$-highest weight elements $b\in B$, then
  $\pr$ is determined on all of $B$.
\end{lemma}
\begin{proof}
  Any element $b'\in B$ is connected to a $\{1,2,\ldots,n-1\}$ highest weight element $b$
  using a sequence $e_{i_1} \cdots e_{i_k}$ with $i_j\in \{1,2,\ldots,n-1\}$. Hence
  $\pr(b')=e_{i_1+1}\cdots e_{i_k+1}(\pr(b))$, which is determined if $\pr(b)$ is known,
  since $i_j+1\in \{2,\ldots,n\}$.
\end{proof}  

\begin{definition}
The \textit{orbit} of $b\in B$ under the promotion operator $\pr$ is the family
\begin{equation*}
  b \prto \pr(b) \prto \pr^2(b) \prto \cdots \prto \pr^n(b) \prto b\,,
\end{equation*}
(or any cyclic shift thereof).
\end{definition}

\begin{lemma}[Two Path Lemma] \label{lemma.twoPaths}
  Suppose $x,y,b\in B$ such that the following conditions hold:
  \begin{enumerate}
    \item The entire orbits of $x$ and $y$ are known;
    \item $b$ is connected to $x$ by a chain of crystal edges, with all edge colors from some set $I_x$;
    \item $b$ is connected to $y$ by a chain of crystal edges, with all edge colors from some set $I_y$;
    \item $I_x \cap I_y = \emptyset$. 
  \end{enumerate}
  Then the entire orbit of $b$ under promotion is determined.
\end{lemma}
\begin{proof}
  By Remark~\ref{remark.shift}, we have $\pr^k \circ e_i = e_{i+k} \circ \pr^k$ (and similarly for
  $f_i$) whenever $i,i + k \neq 0 \pmod{n+1}$. Since by assumption the entire orbit of $x$ is known
  and $b$ is connected to $x$ by a chain consisting of edges from the set $I_x$, all powers
  $\pr^k(b)$ are determined except for $k \in \{n+1-i \}_{i \in I_x}$.  Similarly, the entire orbit
  of $y$ is known and $b$ is connected to $y$ by a chain consisting of edges from the set $I_y$, all
  powers $\pr^k(b)$ are determined except for $k \in \{ n+1-i\}_{i \in I_y}$.  Since $I_x\cap
  I_y=\emptyset$, the entire orbit of $b$ is determined.
\end{proof}

\subsection{Kashiwara's conjecture}
Let $B$ be a $U_q(\mathfrak{g})$-crystal with index set $I$ (for the purpose of this paper it 
suffices to assume that $\mathfrak{g}$ is of type $A_n^{(1)}$, but the statements in this subsection 
hold more generally). 
We denote by $\mathfrak{g}_J$ the subalgebra of $\mathfrak{g}$ restricted to the index set $J\subset I$.
The crystal $B$ is said to be \emph{regular} if, for any $J\subset I$ of finite-dimensional type,
$B$ as a $U_q(\mathfrak{g}_J)$-crystal is isomorphic to a crystal associated with an
integrable $U_q(\mathfrak{g}_J)$-module. Stembridge~\cite{Stembridge.2003} provides a local
characterization of when a $\mathfrak{g}$-crystal is a crystal corresponding to a
$U_q(\mathfrak{g})$-module.

In~\cite{Kashiwara.1994,Akasaka.Kashiwara.1997}, Kashiwara defined the notion of 
extremal weight modules. Here we briefly review the definition of an \emph{extremal weight crystal}
$\tilde{B}(\lambda)$ for $\lambda\in P$. Let $W$ be the Weyl group associated to $\mathfrak{g}$ and 
$s_i$ the simple reflection associated to $\alpha_i$.  Let $B$ be a crystal corresponding to an 
integrable $U_q(\mathfrak{g})$-module. A vector $u_{\lambda}\in B$ of weight $\lambda\in P$ 
is called an \emph{extremal vector} if there exists a family of vectors $\{u_{w\lambda}\}_{w\in W}$ 
satisfying 
\begin{align}
&u_{w\lambda}=u_{\lambda}\text{ for }w=e,\\
&\text{if }\langle \alpha^\vee_i,w\lambda\rangle\ge0,\text{ then }e_iu_{w\lambda}=\fail\text{ and }
f_i^{\langle \alpha^\vee_i,w\lambda\rangle}u_{w\lambda}=u_{s_iw\lambda},\\
&\text{if }\langle \alpha^\vee_i,w\lambda\rangle\le0,\text{ then }f_iu_{w\lambda}=\fail\text{ and }
e_i^{-\langle \alpha^\vee_i,w\lambda\rangle}u_{w\lambda}=u_{s_iw\lambda},
\end{align}
where $\alpha_i^\vee$ are the simple coroots. Then $\tilde{B}(\lambda)$ is an extremal weight crystal 
if it is generated by an extremal weight vector $u_{\lambda}$.

For an affine Kac-Moody algebra $\mathfrak{g}$, let $\delta$ denote the null root in the weight
lattice $P$ and $c$ the canoncial central element. Then define $\Pcl = P/\ZZ\delta$ and
$P^0=\{\lambda \in P \mid \langle c,\lambda \rangle=0 \}$.

\begin{definition}[\cite{Akasaka.Kashiwara.1997}]
  \label{definition.simple}
A finite regular crystal $B$ with weights in $\Pcl^0$ is a \emph{simple crystal} if $B$ satisfies
\begin{enumerate}
\item There exists $\lambda\in \Pcl^0$ such that the weight of any extremal vector of $B$ is
contained in $W_{\mathrm{cl}}\lambda$;
\item The weight space of $B$ of weight $\lambda$ has dimension one.
\end{enumerate}
\end{definition}

\begin{definition}[{Kashiwara~\cite[Section 8]{Kashiwara.2002}}]
\label{definition.good}
  A \emph{`good' crystal} $B$ has the properties that
  \begin{enumerate}
  \item $B$ is the crystal base of a $U_q'(\mathfrak{g})$-module;
  \item $B$ is simple;
  \item Convexity condition: For any $i, j\in I$ and $b\in B$, the function 
  $\varepsilon_i(f_j^k b)$ in $k$ is convex.
  \end{enumerate}
\end{definition}

Note that the third condition of Definition~\ref{definition.good} is only necessary for rank 2 crystals. 
For higher rank crystals this follows from regularity and Stembridge's local characterization of 
crystals~\cite{Stembridge.2003}.

\begin{conjecture}[{Kashiwara~\cite[Introduction]{Kashiwara.2005}}]
\label{conjecture.Kashiwara}
Any `good' finite affine crystal is the tensor product of Kirillov-Reshetikhin crystals.
\end{conjecture}

\section{Promotion}
\label{section.promotion}

In this section we introduce the Sch\"utzenberger operator $\prA$  involving jeu-de-taquin on 
highest weight crystals $B(\lambda)$. This is used to show that promotion operators exist on 
$B(\lambda)$ if and only if $\lambda$ is a rectangle. We then extend the definition of $\prA$ to
tensor products and discuss its relation to connectedness.

\subsection{Existence and uniqueness on $B(\lambda)$} 
\label{subsection.promotion}
Sch\"utzenberger~\cite{Schuetzenberger.1972} defined a weak promotion operator 
$\prA$ on standard tableaux. Here we define the obvious extension~\cite{Shimozono.2002} on 
semi-standard tableaux on the alphabet $\{1,2,\ldots,n+1\}$ using 
jeu-de-taquin~\cite{Schuetzenberger.1977} (see for example also~\cite{Fulton.1997}):
\begin{enumerate}
\item Remove all letters $n+1$ from tableau $t$ (this removes a horizontal strip from $t$);
\item Using jeu-de-taquin, slide the remaining letters into the empty cells (starting from
left to right);
\item Fill the vacated cells with zeroes;
\item Increase each entry by one.
\end{enumerate}
The result is denoted by $\prA(t)$.

\begin{example}
Take $n=3$. Then
\begin{equation*}
         \def\t{\tableau[sbY]}
         \let\c=\vcenter
	t = \c{\t{3,4,4|2,3,3|1,1,2}} \quad \overset{(1)+(2)}{\longrightarrow} \quad
	\c{\t{3,3,3|1,2,2|\bullet,\bullet,1}} \quad \overset{(3)+(4)}{\longrightarrow} \quad
	\c{\t{4,4,4|2,3,3|1,1,2}} = \prA(t).
\end{equation*}
\end{example}

One can consider the reverse operation (which is also sometimes called \emph{demotion}):
\begin{enumerate}
\item Remove all letters $1$ from tableau $t$ (this removes the first
  part of the first row);
\item Using jeu-de-taquin, slide the remaining letters into the empty cells;
\item Fill the vacated cells with $n+2$s;
\item Decrease each entry by one.
\end{enumerate}
The result is denoted by $\prA^{-1}(t)$. We will argue in the proof of the following proposition
why these operations are actually well-defined and inverses of each other.

\begin{proposition}
  \label{proposition.prA} 
  Let $\lambda$ be a partition with at most $n$ parts and let
  $B(\lambda)$ be a type $A_n$ highest weight crystal. Then, $\prA$ is
  the unique weak promotion operator on $B(\lambda)$. Furthermore,
  $\prA$ is a promotion operator if and only if $\lambda$ is a rectangle.
\end{proposition}
Using standardization, the second part of the proposition follows from results of 
Haiman~\cite{Haiman.1992}  who shows that, for standard tableaux on $n+1$ letters, $\prA$
has order $n+1$ if and only if $\lambda$ is a rectangle (and provides a
generalization of this statement for shifted shapes). Shimozono~\cite{Shimozono.2002} 
proves that $\prA$ is the unique promotion operator on $B(s\omega_r)$ of type $A_n$. The
resulting affine crystal is the Kirillov-Reshetikhin crystal $B^{r,s}$
of type $A_n^{(1)}$~\cite{KKMMNN.1992,Shimozono.2002}. 
We could not find the statement of the uniqueness of the weak promotion operator in the 
literature.

For the sake of completeness, we include a complete and elementary proof of 
Proposition~\ref{proposition.prA}; the underlying arguments are similar in spirit to those 
in~\cite{Haiman.1992}, except that we are using crystal operations on semi-standard tableaux 
instead of dual equivalence on standard tableaux. We first recall the following properties of
jeu-de-taquin (see for example~\cite{Kashiwara.1995,LLT.1995,Fulton.1997,
Lascoux.Leclerc.Thibon.1997}).

\begin{remarks}
  \label{recalls.jeu_de_taquin}
  Fix the ordered alphabet $\{1,2,\ldots,n+1\}$.

  (a) Jeu-de-taquin is an operation on skew tableaux which commutes
  with crystal operations.

  (b) Let $\lambda/\mu$ be a skew partition, and $T$ the set of
  semi-standard skew-tableaux of shape $\lambda/\mu$, endowed with its
  usual type $A_n$ crystal structure. Let $f$ be a function which maps
  each skew tableau in $T$ to a semi-standard tableau of partition shape, and which
  commutes with crystal operations. For example, one can take for $f$
  the straightening function which applies jeu-de-taquin to $t\in T$ until it has partition shape.
  Let $C$ be a connected crystal component of $T$. Then, by commutativity with
  crystal operations, there exists a unique partition $\nu$ such that
  $f(C)$ is the full type $A_n$ crystal $B(\nu)$ of tableaux over this
  alphabet. Since $B(\nu)$ has no automorphism, this isomorphism is
  unique, and $f$ has to be straightening using jeu-de-taquin.

  (c) Let $\lambda$ be a rectangle, and $\mu\subset\lambda$.  Consider the
  complement partition $\mu^\complement$ of $\mu$ in the
  rectangle $\lambda$. Then, the type $A_n$ crystal of skew tableaux of
  shape $\lambda/\mu$ is isomorphic to the crystal of tableaux of
  shape $\mu^\complement$; this can be easily seen by rotating each tableau
  $t$ of shape $\mu^\complement$ by $180^\circ$ and mapping each letter $i$
  to $n+2-i$. By uniqueness of the isomorphism, the isomorphism and its inverse are
  both given by applying jeu-de-taquin, either sliding up or down.  In
  particular, jeu-de-taquin down takes any tableau of shape
  $\lambda/\mu$ to a tableau of shape $\mu^\complement$, and vice-versa.
\end{remarks}

\begin{example}
  \label{example.jeu_de_taquin_on_rectangle}
  Let $\lambda:= (6^4)$ and $\mu:=(5,2)$. The complement partition of
  $\mu$ in $\lambda$ is $\mu^\complement=(6,6,4,1)$. We now apply jeu-de-taquin
  up from a tableau of shape $\mu^\complement$, and obtain a skew tableau of
  shape $\lambda/\mu$. Applying jeu-de-taquin down yields back the
  original tableau. As is well-known for jeu-de-taquin, the end result does
  not depend on the order in which the inner corners are filled; here we show one
  intermediate step, after filling successively the three inner corners $(2,4),(5,3)$, and $(6,3)$. 
  The color of the dots at the bottom (resp. at the top) indicates at which step each empty cell has 
  been created by jeu-de-taquin up (resp. down).
  \begin{equation*}
    \let\b\bullet\let\c\circ
    \vcenter{\tab{
      6 ,\c,\c,\b,\b,\b|
      4 ,4 ,5 ,5 ,\c,\c|
      3 ,3 ,3 ,3 ,4 ,6 |
      2 ,2 ,2 ,2 ,3 ,4 
    }}
    \overset{J-D-T}{\longleftrightarrow}
    \vcenter{\tab{
      4 ,6 ,\c,\c,\c,\c|
      3 ,3 ,4 ,5 ,5 ,6 |
      2 ,2 ,3 ,3 ,4 ,4 |
      \b,\b,\b,2 ,2 ,3 
    }}
    \overset{J-D-T}{\longleftrightarrow}
    \vcenter{\tab{
      3 ,4 ,4 ,5 ,6 ,6|
      2 ,3 ,3 ,3 ,4 ,5|
      \b,\b,2 ,2 ,2 ,4|
      \c,\c,\c,\c,\b,3 
    }} \; .
  \end{equation*}
\end{example}

\begin{proof}[Proof of Proposition~\ref{proposition.prA}]
  We first check that $\prA$ is well-defined; the only non trivial
  part is at step $3$ where we must ensure that the previously vacated
  cells form the beginning of the first row. Fix a partition
  $\lambda$, and consider the set $T$ of all tableaux whose $n+1$s are
  in a given horizontal border strip of length $k$. Step (1) puts them
  in bijection with the tableaux of the type $A_{n-1}$ crystal
  $B(\lambda')$ where $\lambda'$ is $\lambda$ with the border strip
  removed. Let $f$ be the function on $B(\lambda')$ which implements
  the jeu-de-taquin step (2) of the definition of $\prA$. Since
  jeu-de-taquin commutes with crystal operations, $B(\lambda')$ is an
  irreducible crystal, and since crystal operations preserve shape,
  all tableaux in $f(B(\lambda'))$ have the same skew-shape
  $\lambda/\mu$. Considering $f(t)$ where $t$ is the anti-Yamanouchi
  tableau of shape $\lambda'$ shows that $\mu=(1^k)$ as desired
  because the jeu-de-taquin slides follow successive hooks (the \emph{anti-Yamanouchi}
  tableau of shape $\lambda'=(\lambda_1',\ldots,\lambda_m')$ is the unique tableau
  of shape $\lambda'$ which contains $\lambda'_i$ entries $m+1-i$).  
  For example:
  \begin{equation*}
    \let\b=\bullet
    \def\r{\color{DarkRed}}
    \def\Tscale{.65}
    \scriptsize
    \vcenter{\tab{
        4 ,  |
        3 ,4 ,4 ,  ,  |
        2 ,3 ,3 ,4 ,4 ,  |
        1 ,2 ,2 ,3 ,3 ,4 ,4}}
    \rightarrow
    \vcenter{\tab{
        4 ,\r4  |
        3 ,\r3 ,4 ,  ,  |
        2 ,\r2 ,3 ,4 ,4 ,  |
        \r\b,\r1 ,2 ,3 ,3 ,4 ,4}}
    \rightarrow
    \vcenter{\tab{
        4 ,4  |
        3 ,3 ,4  ,\r4,  |
        2 ,2 ,3  ,\r3,4 ,  |
        \b,\r\b,\r1,\r2,3 ,4 ,4}}
    \rightarrow
    \vcenter{\tab{
        4 ,4  |
        3 ,3 ,4 ,4  ,\r4|
        2 ,2 ,3 ,3  ,\r3,  |
        \b,\b,\r\b,\r1,\r2,4 ,4}}
    \rightarrow
    \vcenter{\tab{
        4 ,4  |
        3 ,3 ,4 ,4  ,4  |
        2 ,2 ,3 ,3  ,3  ,\r4|
        \b,\b,\b,\r\b ,\r1,\r2,4}}\; .
  \end{equation*}
  Note further that applying down jeu-de-taquin to $f(t)$ reverses the
  process, and yields back $t$. It follows that $\prA^{-1}$ as
  described above is indeed a left inverse and therefore an inverse
  for $\prA$. Finally, $\prA$ satisfies conditions~\eqref{pt.content}
  and~\eqref{pt.e} of Definition~\ref{definition.promotion} by construction, so it 
  is a weak promotion operator.

  We now prove that a weak promotion operator $\pr$ on $B(\lambda)$ is
  necessarily $\prA$.  Consider the action of $\pr^{-1}$ on a tableau
  $t$. By condition~\eqref{pt.content} of Definition~\ref{definition.promotion}, it has to strip away 
  the $1$s, subtract one from each remaining letter, transform the result into
  a semi-standard tableau of some shape $\mu'(t)\subset \lambda$, and
  complete with $n+1$s. Let $B'$ be the set of all skew-tableaux in
  $B(\lambda)$ after striping and subtraction, endowed with the
  $A_{n-1}$ crystal structure induced by the $\{2,\dots,n\}$ crystal
  structure of $B(\lambda)$. Write $f^{-1}$ for the function which
  reorganizes the letters. By condition~\eqref{pt.e} of Definition~\ref{definition.promotion}, 
  $f^{-1}$ is an $A_{n-1}$-crystal morphism, so by Remark~\ref{recalls.jeu_de_taquin}
  (b) it has to be jeu-de-taquin.  Therefore $\pr^{-1}=\prA^{-1}$, or
  equivalently $\pr=\prA$.

  It remains to prove that $\prA$ is a promotion operator if and only
  if $\lambda$ is a rectangle.

  Assume first that $\lambda$ is a rectangle. By
  Remark~\ref{recalls.jeu_de_taquin} (c), for each $k$, jeu-de-taquin
  down provides a suitable bijection $f^{-1}$ from skew tableaux of
  shape $\lambda/(k)$ and tableaux of shape $(k)^\complement$. The
  inverse bijection $f$ is jeu-de-taquin up. We show $\prA^2 \circ e_n
  = e_1 \circ \prA^2$, which by Remark~\ref{remark.shift2} finishes the
  proof that $\prA$ is a promotion operator.
  Let $t$ be a semi-standard tableau, $l_1$, $l_2$, and $l_3$ be
  respectively the number of bracketed pairs $(n+1,n)$, of unbracketed
  $n+1$s, and unbracketed $n$s. Then, from Remark~\ref{recalls.jeu_de_taquin} (c) one 
  can further deduce that in $\prA^2(t)$ there are $l_1$ bracketed pairs $(2,1)$, $l_2$
  unbracketed $2$s, and $l_3$ unbracketed $1$s. We revisit
  Example~\ref{example.jeu_de_taquin_on_rectangle} in this context. We
  have $l_1=2$, $l_2=1$, and $l_3=2$;
  \let\b\bullet\let\c\circ
  due to label shifts, we have on the left $\b=n+1$ and $\c=n$, in the
  middle $\b=1$ and $\c=n+1$, and on the right $\b=2$ and $\c=1$:
  \begin{equation*}
    \vcenter{\tab{
      5 ,\c,\c,\b,\b,\b|
      3 ,3 ,4 ,4 ,\c,\c|
      2 ,2 ,2 ,2 ,3 ,5 |
      1 ,1 ,1 ,1 ,2 ,3 
    }}
    \overset{\prA}{\longrightarrow}
    \vcenter{\tab{
      4 ,6 ,\c,\c,\c,\c|
      3 ,3 ,4 ,5 ,5 ,6 |
      2 ,2 ,3 ,3 ,4 ,4 |
      \b,\b,\b,2 ,2 ,3 
    }}
    \overset{\prA}{\longrightarrow}
    \vcenter{\tab{
      4 ,5 ,5 ,6 ,7 ,7|
      3 ,4 ,4 ,4 ,5 ,6|
      \b,\b,3 ,3 ,3 ,5|
      \c,\c,\c,\c,\b,3 
    }} \; .
  \end{equation*}
  It follows in particular that $e_1$ applies to $\prA^2(t)$ if and
  only if $l_2>0$ if and only if $e_n$ applies to $t$; furthermore
  both the action of $e_1$ and $e_n$ decrease $l_2$ by one and
  increase $l_3$ by one. This does not change $\mu=(l_1,l_1+l_2+l_3)$,
  and therefore the jeu-de-taquin action on the rest of the
  tableaux. Therefore $\prA^2(e_n(t)) = e_1(\prA^2(t))$, as desired.
  
  To conclude, let us assume that $\lambda$ is not a rectangle. We
  show that $\prA^2 \circ e_n \ne e_1 \circ \prA^2$, which by
  Remark~\ref{remark.shift2} implies that $\prA$ cannot be a promotion
  operator. The prototypical example is $n=2$ and $\lambda=(2,1)$,
  where the following diagram does not commute:
  \begin{equation*}
    \vcenter{
      \begin{tikzpicture}
        \let\u=\underline
        \let\c=\vcenter
        \def\t{\tableau[sbY]}
        \matrix (bla) [matrix of math nodes, column sep={3cm,between origins}, row sep = 1cm]
        {
          {\c{\t{3|\u1,3}}} & {\c{\t{\u2|1,1}}} & {\c{\t{\u3|2,2}}}\\
          {\c{\t{2|\u1,3}}} & {\c{\t{3|1,\u2}}} & {\c{\t{2|1,\u3}} \ne \c{\t{\u3|1,2}} \color{white}{\ne \t{2,2}}}\\
        };
        \draw[->] (bla-1-1) to node [auto] {$\prA$} (bla-1-2);
        \draw[->] (bla-2-1) to node [auto] {$\prA$} (bla-2-2);
        \draw[->] (bla-1-2) to node [auto] {$\prA$} (bla-1-3);
        \draw[->] (bla-2-2) to node [auto] {$\prA$} (bla-2-3);
        
        \draw[->] (bla-1-1) to node [auto] {$e_n$} (bla-2-1);
        \draw[->] (bla-1-3) to node [auto] {$e_1$} (bla-2-3);
      \end{tikzpicture}}
  \end{equation*}
  Interpretation: the underlined cell is the unique cell containing a
  $1$ (resp. a $2$, resp. a $3$) on the left hand side (resp. middle,
  resp. right hand side), and we can track how it moves under
  promotion. Note that the value in the cell is such that promotion
  will always move the cell weakly up or to the right, and neither $e_1$ nor 
  $e_n$ affects it. At the first promotion step, depending on
  whether we apply $e_n$ or not, the cell moves to the right, or
  up. But then, due to the inner corner of the partition it cannot
  switch to the other side, and therefore the diagram cannot close.

  The same phenomenon occurs for any shape having (at least one)
  inner corner. Consider the uppermost inner corner, and construct the
  tableau:
  \begin{equation}
    \def\cell#1#2#3{\multicolumn{1}{#1@{\hspace{.6ex}}c@{\hspace{.6ex}}#2}{\raisebox{-.3ex}{$#3$}}}
    \def\r {\cell{}|}
    \def\lr{\cell||}
    \def\mycdots{\qquad\cdots\qquad}
    \begin{array}{|cccccccccc}
      \cline{1-4}
      \lr{n-1}    & \mycdots    & \lr{n-1}    & \lr{n+1}\\
      \cline{1-4}
      \lr{n-2}    & \mycdots    & \lr{n-2}    & \lr{n-1}\\
      \cline{1-4}
      \lr{\vdots} & \lr{\vdots} & \lr{\vdots} & \lr{\vdots}\\
      \cline{1-4}
      \lr{n-k}    & \mycdots    & \lr{n-k}    & \lr{n-k+1}\\
      \cline{1-6}
                  &             &             & \lr{\underline{n-k}}   & \lr{n+1} & \mycdots\\
      \cline{4-6}
                  &             &             &            &          & \mycdots\\
                  & <n-k        &             &            &          & \mycdots\\
                  &             &             &            &          & \mycdots\\
      \cline{1-6}
    \end{array}
  \end{equation}
  We assume that this tableau does not contain any letter $n$ so that 
  $e_n$ applies to it and transforms the $n+1$ in the top row into an $n$.
  Let $j$ be the width of the upper rectangle and assume that the tableau
  does not contain any further letters $n-k$ (only the $j$ copies in the first $j$ columns).
  Applying $\prA$ without application of $e_n$, promotion slides the
  underlined $n-k$ up, and even after an additional application of $\prA$ all the $j$ cells 
  containing $n-k$ in the original tableau, are in the upper rectangle. First applying $e_n$ and 
  then $\prA$ has the effect of sliding the cell containing the underlined $n-k$
  to the right; this cell cannot come back in the upper rectangle with another application
  of $\prA$. Hence $\prA^2\circ e_n \neq e_1 \circ \prA^2$.
\end{proof}

\subsection{Promotion on tensor products}
Now take $B:=B(s_1\omega_{r_1}) \otimes \cdots \otimes B(s_\ell \omega_{r_\ell})$
a tensor product of $\ell$ classical highest weight crystals labeled by rectangles.
For $b_1\otimes \cdots\otimes b_\ell \in B$, define $\prA: B\to B$ by 
\begin{equation} \label{eq:prAOnTensor}
	\prA(b_1\otimes \cdots \otimes b_\ell) = \prA(b_1) \otimes \cdots \otimes \prA(b_\ell).
\end{equation}

\begin{lemma} \label{lemma.existence}
  $\prA$ on $B=B(s_1\omega_{r_1}) \otimes \cdots \otimes B(s_\ell \omega_{r_\ell})$
  is a connected promotion operator.
\end{lemma}
\begin{proof}
  Since $\prA$ on each tensor factor $B(s_i\omega_i)$ satisfies conditions~\eqref{pt.content}
  and~\eqref{pt.order} of Definition~\ref{definition.promotion}, $\prA$ on $B$ also satisfies 
  conditions~\eqref{pt.content} and~\eqref{pt.order}. Since $\prA$ on each tensor factor 
  $B(s_i\omega_i)$ satisfies condition~\eqref{pt.e} and the bracketing is well-behaved with 
  respect to acting on each tensor factor, we also have condition~\eqref{pt.e} for $\prA$ on $B$. The 
  affine crystal resulting from $\prA$ on $B(s\omega_r)$ is the Kirillov-Reshetikhin crystal $B^{r,s}$ of 
  type $A_n^{(1)}$~\cite{KKMMNN.1992,Shimozono.2002}. Since $B^{r,s}$ is simple,
  the affine crystal resulting from $\prA$ on $B$ is connected 
  by~\cite[Lemmas 4.9 and 4.10]{Kashiwara.2002}.
\end{proof}

Lemma~\ref{lemma.existence} shows that a promotion operator with the properties of
Definition~\ref{definition.promotion} exists on $B=B(s_1\omega_{r_1}) \otimes \cdots \otimes 
B(s_\ell \omega_{r_\ell})$. Theorem~\ref{theorem.main} states that for $\ell=2$ this is 
the only connected promotion operator.

\section{Inductive proof of Theorem~\ref{theorem.main}}
\label{section.proof}
In this section we provide the proof of Theorem~\ref{theorem.main}. Throughout this section 
$B:=B(s'\omega_{r'}) \otimes B(s\omega_r)$. For $n<\max(r,r')$ this crystal is 
either nonexistent or trivial.

\subsection{Outline of the proof}
Aside from distinguishing the cases where $(s',r') = (s,r)$, our proof does not depend in a material
way on the values of $s$ and $s'$.  The basic tool in our proof is an induction which
allows us to relate the cases described by the triple $(r',r,n)$ to those described by
$(r'-1,r-1,n-1)$, provided that 
\begin{enumerate}
  \item $n \le r' + r -1$ and
  \item we do not have $r' = r = 1$. 
\end{enumerate}
As follows from Lemma~\ref{lemma.duality}, any crystal which does not satisfy these hypotheses
is isomorphic to one which does, with the exception of the case where $r' = r = n = 1$. (This case
does not satisfy the result of the Theorem as was discussed in Remark~\ref{remark.counterexample}).
The general idea for the proof is to use repeated applications
of induction and duality to successively reduce the rank of the crystal.  Note that both techniques
preserve the fact that the rank is greater or equal to the maximum of the heights of the two rectangles 
$r$ and $r'$. We take as base cases those crystals where either $r$ or $r'$ is equal to zero.  
In these cases, we have only a single tensor factor and the statement of Theorem~\ref{theorem.main}
was shown by Shimozono~\cite{Shimozono.2002}.

This approach, however, does not cover those cases which inductively reduce to the case $(1,1,1)$.
The only case which \emph{directly} reduces to $(1,1,1)$ is $(2,2,2)$.  By duality, the
case $(2,2,2)$ is equivalent to the case $(1,1,2)$.  We prove this case directly, as a separate base
case, and thus complete the proof.

The proof is laid out as follows.  In Section~\ref{subsection.rowRow}, we discuss the base case of
the $A_2$ crystals with $r = r' = 1$.  In Section~\ref{subsection.induction}, we present the basic
lemma (Lemma~\ref{Prop:Induction}) for our inductive arguments. In
Section~\ref{subsection.rectangle}, we show how to apply the induction in the case where $r' \ge
r$ and $r' > 1$ for different tensor factors. Note that by~\eqref{eq:unique_iso} we can always assume 
that $r'\ge r$. In Section~\ref{subsection.equal_factors} we treat the case of equal tensor factors.

\subsection{Row tensor row case, $n=2$} 
\label{subsection.rowRow}

In this subsection we prove Theorem~\ref{theorem.main} for the row
tensor row case with $n=2$. In this case, the isomorphism of
Equation~\eqref{equation.isomorphism} becomes:
\begin{equation}
  \Psi: B(s\omega_1)\otimes B(s\omega_1)
  \hookrightarrow\!\!\!\!\!\rightarrow
  B((s-1)\omega_1) \otimes B((s+1)\omega_1)\oplus B(s\omega_2)\,.
\end{equation}

\begin{proposition} 
  \label{proposition.rowxrow}
  Let $B := B(s'\omega_1) \otimes B(s\omega_1)$ be the tensor product
  of two single row classical highest weight crystals of type $A_2$ with $s,s'\ge 1$.
  If $s\ne s'$, there is a unique promotion operator
  $\pr=\prA$ which is connected. If $s=s'$, there are two
  promotion operators: $\prA$ which is connected, and $\prA':=\Psi^{-1}\circ \prA \circ \Psi$, 
  induced by the canonical promotions on the classical crystals
  $B((s-1)\omega_1) \otimes B((s+1)\omega_1)$ and $B(s\omega_2)$,
  which is disconnected.
\end{proposition}

We may assume without loss of generality that $s'\le s$. After a
preliminary Lemma~\ref{lemma.content}, we show that if the
$\pr$-orbits coincide with the $\prA$-orbits on the inversionless
component, then $\pr=\prA$ (Proposition~\ref{proposition.inversionlessToAnywhere}). 
Here the inversionless component is the component $B((s+s')\omega_1)$ in 
the decomposition~\eqref{eq:classicalDecomposition} of $B$.
Then, we proceed with the analysis of $\pr$-orbits on the inversionless
component (Lemma~\ref{lemma.inversionless}). When $s'<s$, there is a
single possibility which implies $\pr=\prA$. When $s'=s$, there are two
possibilities, and we argue that one implies $\pr=\prA$, while the other
implies $\pr=\prA'$ via the isomorphism $\Psi$.

\begin{lemma}[Content Lemma] \label{lemma.content}
  If $v' \otimes v\in B$ does not contain any $3$s, then $\pr(v' \otimes v)=\prA(v')\otimes \prA(v)$.  
\end{lemma}
\begin{proof}
  By assumption $w=v'\otimes v$ contains only the letters 1 and 2. The only $1$-bracketing
  that can be achieved is by 2s in the left tensor factor that bracket with 1s in the right tensor
  factor.  Hence knowing $\varphi_1(w)$ and $\varepsilon_1(w)$ determines $w$ completely. 
  Since $\pr$ rotates content, $\pr(w)$ contains only 2s and 3s. Since furthermore
  $\varphi_2(\pr(w))=\varphi_1(w)$ and $\varepsilon_2(\pr(w))=\varepsilon_1(w)$, this completely
  determines $\pr(w)$.  Since $\prA$ is a valid promotion operator by Lemma~\ref{lemma.existence}, 
  $\pr$ must agree with $\prA$ on these elements.
\end{proof}

\begin{proposition}
  \label{proposition.inversionlessToAnywhere}
  Let $B := B(s'\omega_1) \otimes B(s\omega_1)$, and $\pr$ be a
  promotion on a classical type $A_n$ crystal $C := B\oplus B'$ of which
  $B$ is a direct summand (typically $C:=B$). Assume that the orbits
  under promotion on the inversionless component of $B$ coincide with
  those for the canonical promotion $\prA$ of $B$. Then $\pr$
  coincides with $\prA$ on $B$.
\end{proposition}
We start with the elements with only one letter in some tensor factor.
\begin{lemma}
  \label{lemma.oneLetter} 
  Under the hypothesis of
  Proposition~\ref{proposition.inversionlessToAnywhere}, the
  $\pr$-orbit of an element $v' \otimes v$ is its $\prA$-orbit
  whenever either $v'$ or $v$ contains a single letter.
\end{lemma}
\begin{proof}
  Assume that $v'=k^{s'}$ (resp. $v=k^s$). Then, $v' \otimes v$ is in
  the $\prA$-orbit of the inversionless element $1^{s'}\otimes
  \prA^{1-k}(v)$ (resp. of $\prA^{3-k}(v')\otimes 3^{s}$) which by hypothesis
  is also its $\pr$-orbit.
\end{proof}
Next come elements with exactly two letters in each tensor factor.
\begin{lemma}
  \label{lemma.twoLetter} 
  Under the hypothesis of
  Proposition~\ref{proposition.inversionlessToAnywhere}, the
  $\pr$-orbit of an element $w:=v' \otimes v$ is its $\prA$-orbit
  whenever $v'$ and $v$ each contain precisely two distinct letters.
\end{lemma}
\begin{proof}
  By Lemma~\ref{lemma.oneLetter}, it remains to consider the cases when both $v'$ and $v$
  contain two letters.
  \begin{enumerate}
  \item If $w=1^a2^b \otimes 2^c3^d$ it is inversionless and we are done by hypothesis.
    The $\prA$-orbit includes the elements $2^a3^b \otimes 1^d 3^c$ and $1^b3^a \otimes 1^c2^d$. 
  \item Assume $w = 2^a 3^b \otimes 1^d 2^c$.  Applying $f_2$ a
    sufficient number of times gives $3^{a + b} \otimes 1^d 2^{c_1}
    3^{c_2}$. If we instead apply $e_1$ a sufficient number of times
    to $w$, we get the elements $1^{a_1} 2^{a_2} 3^b \otimes 1^{d +
      c}$. In both cases Lemma~\ref{lemma.oneLetter} applies, and by
    the Two Path Lemma~\ref{lemma.twoPaths}, the $\pr$-orbit of $w$ is its
    $\prA$-orbit.
  \item The orbits of the elements considered previously include the elements $1^b 3^a \otimes 2^d
    3^c$ and $1^a2^b \otimes 1^c3^d$.
  \item Assume $w= 1^b 3^a \otimes 1^d 3^c$.  Applying $e_2^a$ yields
    $1^b 2^a \otimes 1^d 3^c$, and applying $f_1^d$ yields $1^b 3^a
    \otimes 2^d 3^c$. Both elements have already been treated, and by
    the Two Path Lemma~\ref{lemma.twoPaths}, the $\pr$-orbit of $w$ is its $\prA$-orbit.
  \item The orbits of the elements considered previously include $w= 1^a2^b \otimes 1^c2^d$ and $w=
    2^a3^b \otimes 2^c3^d$. Hence all cases are covered.
  \end{enumerate}
\end{proof}

We are now in the position to prove Proposition~\ref{proposition.inversionlessToAnywhere}.
\begin{proof}[Proof of Proposition~\ref{proposition.inversionlessToAnywhere}]
  By the Highest Weight Lemma~\ref{lemma.highestWeight}, we only need
  to determine promotion of each $\{1\}$-highest weight element. They are of the form 
  $1^a 2^b 3^c \otimes 1^d 3^e$ with $b \le d$.
  We claim that its promotion orbit is given as follows:
  \begin{equation*}
    w_0 = 1^a 2^b 3^c \otimes 1^d 3^e \stackrel{(1)}{\longrightarrow} w_1 = 1^c 2^a 3^b \otimes 1^e 2^d 
    \stackrel{(2)}{\longrightarrow} w_2 = 1^b 2^c 3^a \otimes 2^e 3^d 
    \stackrel{(3)}{\longrightarrow} w_0 \; .
  \end{equation*}
  Applying $e_1$ a sufficient number of times to $w_1$ yields a
  word whose second factor contains a single letter. Using
  Lemma~\ref{lemma.oneLetter}, we deduce that
  $\pr(w_1)=w_2=\prA(w_1)$ as claimed (arrow (2)).
  Applying $f_1^b$ to $w_2$ yields $2^{b+c} 3^a \otimes 2^c 3^d$ whose
  $\pr$-orbit is its $\prA$-orbit by Lemma~\ref{lemma.twoLetter}. Therefore
  $\pr(w_2)=w_0=\prA(w_2)$ as claimed (arrow (3)).
  Arrow (1) follows from $\pr^3=\id$.
\end{proof}

We now turn to the analysis of promotion orbits on the inversionless
component (see Figure~\ref{fig:inversionlessComponent}).
\begin{figure}
  \begin{bigcenter}
    \newcommand{\myinput}[1]{\scalebox{0.4}{\input{#1}}}
    \begin{tabular}{l@{\qquad}l}
      \scalebox{0.4}{\input{figure-row3xrow2-inversionLess}} &
      \scalebox{0.4}{\input{figure-row3xrow3-inversionLess}}\\\\
      \scalebox{0.4}{\input{figure-row4xrow2-inversionLess}} &
      \scalebox{0.4}{\input{figure-row5xrow1-inversionLess}}
    \end{tabular}
  \end{bigcenter}
  \caption{The inversionless classical component of the $B(s'\omega_1)
    \otimes B(s\omega_1)$ $A_2$-crystal. The $\pr$-orbit is uniquely
    determined by content for the grayed elements, by bracketing for
    the framed elements, and then by the Two Path Lemma for the
    lightly framed elements. This covers all elements except when
    $s=s'$ (upper right), where there are potentially two choices for
    $1^s\otimes 3^s$, and consequently for any element in its row or
    column.}
  \label{fig:inversionlessComponent}
\end{figure}
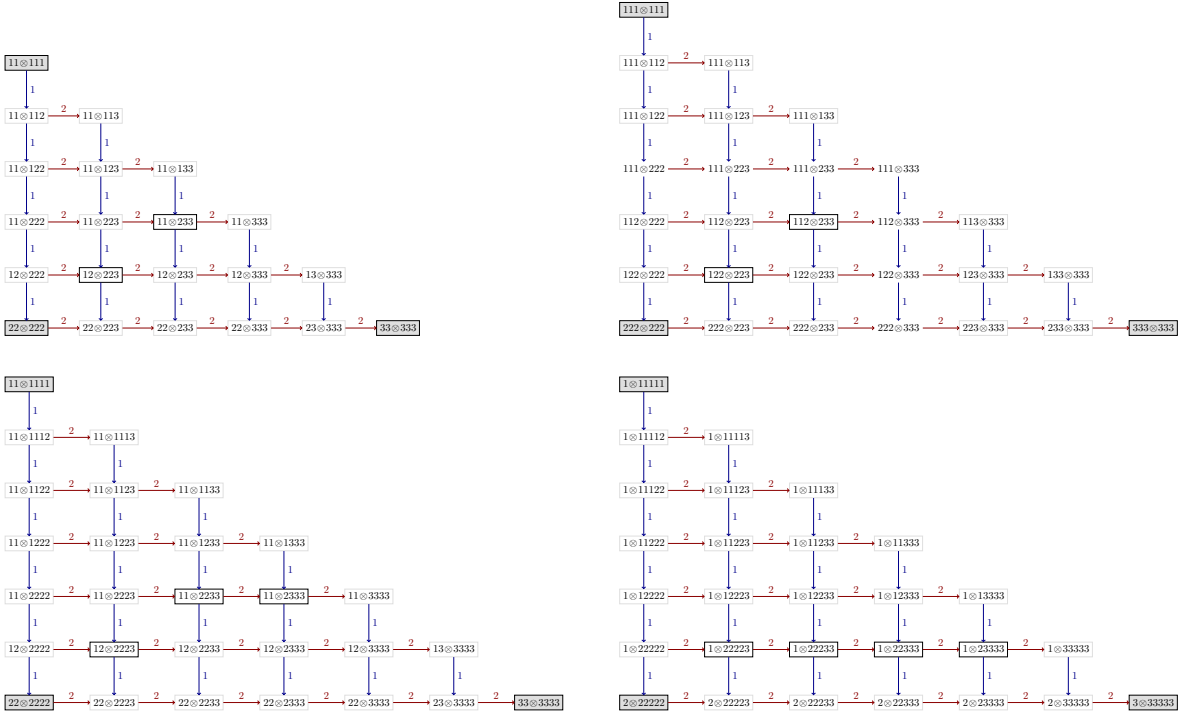
\begin{lemma} 
  \label{lemma.inversionless}
  When $s'\ne s$, the $\pr$-orbit of every element in the inversionless 
  component agrees with $\prA$.  When $s' = s$, there are precisely two 
  cases; either $\pr$ agrees with $\prA$ on the orbit of every element in 
  this component, or $\pr$ agrees with $\prA'$ on the orbit of every
  element in this component.
\end{lemma}

\begin{proof}
  Draw the crystal graph for the inversionless component with
  $1^{s'}\otimes 1^s$ at the top, $1$-arrows going down and $2$-arrows
  going right (see Figure~\ref{fig:inversionlessComponent}).  When there is no ambiguity,
  we drop the $\otimes$ sign and consider elements in $B$ as words. The orbits of the elements 
  $w$ in the inversionless component are considered in the following order:
  \begin{enumerate}
  \item Corners: $w\in \{1^{s'}\otimes 1^s,2^{s'} \otimes 2^s, 3^{s'} \otimes 3^s\}$.
  \item Diagonal: $w:=1^a 2^{s'-a} \otimes 2^{s-a} 3^a$ with $1\le a<s'$.
  \item Middle row: $w:=1^{s'} \otimes 2^{s-a} 3^a$ with $s' \le a < s$.
  \item Lower leftmost column: $w:=1^a 2^{s'-a} \otimes 2^s$ with $1\le a \le s'$ and $a<s$.
  \item Left of lower row: $w:=2^{s'} \otimes 2^a 3^{s-a}$ with $1\leq a<s$.
  \item Rest of leftmost column and lower row: $w:=1^a 2^b$ or $w:=2^a
    3^b$, except when $a=b=s=s'$.
  \item General elements: $w:=1^a 2^b 3^c$ not in any of the other cases.
  \item Row and column of $1^s\otimes3^s$ when $s=s'$.
  \end{enumerate}

  (1): By content, $1^{s'} \otimes 1^s \prto 2^{s'} \otimes 2^s \prto
  3^{s'} \otimes 3^s \prto 1^{s'} \otimes 1^s$, which agrees with the $\prA$-orbit.

  (2): The orbit of $w:=1^a 2^{s'-a} \otimes 2^{s-a} 3^a$ for $1\le
  a<s'$ is forced by bracketing arguments. Recall that $b_i(w)$ denotes the number of 
  $()$ brackets in the construction of $f_i$ and $e_i$ on $w$. Start with the element
  $w_1 := 2^a 3^{s' - a} \otimes 1^a 3^{s-a}$. Note that
  $b_1(w_1) = a$.  Thus $b_2(\pr(w_1)) = a$.  This
  implies that in $w_2 := \pr(w_1)$ all $3$s must be
  in the left tensor factor and all $2$s must be in the right tensor
  factor.  This forces $w_2 = 1^{s'-a} 3^a \otimes 1^{s-a}
  2^a$.  Now we have $b_1(w_2) = 0$. Thus if we define
  $w_0:= \pr(w_2) = \pr^{-1}(w_1)$, we must have $b_2(w_0) = 0$.
  However, we also have that $b_2(w_1) = 0$.  Hence we have $b_1(w_0)
  = 0$. These facts imply that
  in $w_0$ all $1$s precede all $2$s and all $2$s precede all $3$s.
  Therefore $w_0 = 1^a 2^{s' - a} \otimes 2^{s - a}
  3^a = w$, and the $\pr$-orbit of $w$ is its $\prA$-orbit $w \prto w_1
  \prto w_2 \prto w$.
    
  (3): The argument for $w:=1^{s'} \otimes 2^{s-a} 3^a$ with
  $s' \le a < s$ is very similar to (2).  Start with $w_1:=2^{s'} \otimes 1^a
  3^{s-a}$. Since $b_1(w_1) = s'$, we must have $b_2(w_2) = s'$ for 
  $w_2:=\pr(w_1)$. This
  forces the first tensor factor of $w_2$ to be $3^{s'}$, which
  completely fixes $w_2=3^{s'} \otimes 1^{s-a} 2^a$. We have $b_1(w_2) = 0$, and so $b_2(w_0)
  = 0$ where $w_0 := \pr(w_2) = \pr^{-1}(w_1)$. However, we also have $b_2(w_1) = 0$, and so
  $b_1(w_0) = 0$.  Thus, as above, $w_0$ must have no inversions, and we have
  $w_0 = 1^{s'} \otimes 2^{s-a} 3^a$ and the $\pr$-orbit of $w$ is its $\prA$-orbit.
    
  (4): Applying $f_1^a$ to $w:=1^a 2^{s'-a} \otimes 2^s$ gives the
  corner element $2^{s'} \otimes 2^s$. Applying $f_2^a$ to $w$ yields
  the diagonal element $1^a 2^{s'-a} \otimes 2^{s-a} 3^a$ (or middle
  row element for $a=s'<s$). Hence by the Two Path
  Lemma~\ref{lemma.twoPaths} the $\pr$-orbit of $w$ is its $\prA$-orbit.
    
  (5): Let $w:=2^{s'} \otimes 2^{s-a} 3^a$ for $1 \le a < s$,
  Then $e_2^a(w)=2^{s'} \otimes 2^s$ whose orbit is known by the corner
  case, and $e_1^{s'}$ or $e_1^a$ applied to $w$ yields a diagonal or row
  element. Hence again by the Two Path Lemma~\ref{lemma.twoPaths} the
  $\pr$-orbit of $w$ is its $\prA$-orbit.
        
  (6): By Lemma~\ref{lemma.content}, the element $w:=1^a 2^b$ is
  mapped to $2^a 3^b$ under promotion. From Step (4) and (5) either
  the $\pr$-orbit of $1^a 2^b$ or of $2^a 3^b$ is its $\prA$-orbit,
  except when $a=b=s=s'$.
   
  (7): For a general element $w:=1^a 2^b 3^c$, $f_1^a(w)$ yields the
  element $2^{a+b} 3^c$ of the lowest row, and $e_2^c(w)$ yields the
  element $1^a 2^{b+c}$ of the leftmost column edge. By (4), (5), (6)
  and the Two Path Lemma~\ref{lemma.twoPaths} the $\pr$-orbit of $w$ is
  its $\prA$-orbit, except when $s=s'$ and $w$ is in the row or column
  of $1^s\otimes3^s$.

  (8): Assume $s=s'$, and consider $\pr(1^{s-1} 2 \otimes 3^s)$. Using
  (7) for the $\pr$-orbit of $1^{s-1}3\otimes 3^s$ and $\pr^3=\id$, we
  have:
    \begin{equation}
      \begin{aligned}
        \pr(          1^{s-1}2\otimes 3^s)   &= 
        \pr(e_2 (     1^{s-1}3\otimes 3^s))  &&=
        \pr(e_2 (\pr (23^{s-1}\otimes 2^s))) \\&=
        \pr^2   (e_1 (23^{s-1}\otimes 2^s))  &&=
        \pr^{-1}(e_1 (23^{s-1}\otimes 2^s))  \\&=
        \pr^{-1}(     13^{s-1}\otimes 2^s)\,.
      \end{aligned}
    \end{equation}
    By the bracketing structure between the $2$s and $3$s in $1
    3^{s-1} \otimes 2^s$, there are only two possible
    choices for $\pr^{-1}(13^{s-1} \otimes 2^s)$, as desired:
    \begin{equation}
      2^{s-1} 3 \otimes 1^s \text{ and } 1 2^{s-1} \otimes 1^{s-1} 3\,.
    \end{equation}
    By the Two Path Lemma~\ref{lemma.twoPaths}, each of those choices
    completely determines the $\pr$-orbits on the row and column of
    $1^s\otimes3^s$. The first choice is clearly compatible with $\prA$. 
    The second choice is
    compatible with $\prA'$: 
    \begin{displaymath}
      12^{s-1} \otimes 1^{s-1} 3 \quad \overset{\Psi}{\longrightarrow} \quad 2^{s-1} \otimes 1^s 3
      \quad \overset{\prA}{\longrightarrow} \quad 3^{s-1} \otimes 1 2^s \quad 
      \overset{\Psi^{-1}}{\longrightarrow} \quad 13^{s-1} \otimes 2^s.
    \end{displaymath}
\end{proof}

\begin{proof}[Proof of Proposition~\ref{proposition.rowxrow}]
  If $s\neq s'$, the statement of the Proposition follows from Lemma~\ref{lemma.inversionless} 
  and Proposition~\ref{proposition.inversionlessToAnywhere}.
  It remains to settle the case $s=s'$.

  Case 1: $\pr$ coincides with $\prA$ on the inversionless
  component. Then, by
  Proposition~\ref{proposition.inversionlessToAnywhere}: $\pr=\prA$.

  Case 2: $\pr$ coincides with $\prA'$ on the inversionless
  component. Note that $\Psi$ maps this component to the inversionless
  component of $B((s-1)\omega_1) \otimes B((s+1)\omega_1)$. Therefore,
  applying Proposition~\ref{proposition.inversionlessToAnywhere} to
  $B((s-1)\omega_1) \otimes B((s+1)\omega_1)$ yields that $\pr$
  coincides with $\prA'$ on $\Psi^{-1}(B((s-1)\omega_1) \otimes
  B((s+1)\omega_1))$. Then, $\pr$ stabilizes both
  $\Psi^{-1}(B((s-1)\omega_1) \otimes B((s+1)\omega_1))$ and
  $\Psi^{-1}(B(s\omega_2))$; since the latter piece admits a unique
  promotion, $\pr=\prA'$.

  Since $\prA\ne\prA'$, the two cases above are necessarily exclusive,
  and by Lemma~\ref{lemma.inversionless} they cover all the choices
  for $\pr$ on the inversionless component.
\end{proof}

\subsection{Induction} \label{subsection.induction}
The remainder of our proof uses induction to relate any promotion operator on 
$B:=B(s'\omega_{r'}) \otimes B(s \omega_r)$ of type $A_n$ to a promotion operator on 
$D:=B(s' \omega_{r'-1}) \otimes B(s \omega_{r-1})$ of type $A_{n-1}$.
This is done in Proposition~\ref{Prop:Induction}. Before we can state and prove this
proposition we first need some more notation and a preliminary lemma. 

Let $C_i(B)$ be the subgraph of $B$ consisting of the vertices with $s+s'$ copies of the letter 
$i$, along with the arrows $e_j,f_j$ with $j\in \{1,2,\ldots,n\} \setminus \{i-1,i\}$. Recall the
promotion operator $\prA$ of Section~\ref{subsection.promotion}. In addition, if $(s,r)=(s',r')$,
we can define the promotion operator $\prA'=\Psi^{-1} \circ \prA \circ \Psi$. If we want to
emphasize that these promotion operators act on the crystal $B$ (respectively $D$), we write 
$\prA_B$ and $\prA'_B$ (respectively $\prA_D$ and $\prA'_D$).

\begin{lemma} \label{lemma.prAandprA'}
  If $(s,r)=(s',r')$, there exist at least two promotion operators $\prA_B$ and $\prA'_B$
  on $B$. Furthermore, they are distinct when restricted to maps on $C_1(B) \to C_2(B)$ or
  $C_n(B) \to C_{n+1}(B)$.
\end{lemma}

\begin{proof}
  Set $C_i:=C_i(B)$.
  We note first that for content reasons, any promotion operator must bijectively map $C_i \to
  C_{i+1}$ (and $C_{n+1} \to C_1$).  By Section~\ref{subsection.promotion} we know that there are at
  least two promotion operators $\prA_B$ and $\prA_B'$ defined on $B$. It remains to show they
  differ when restricted to $C_1$ and $C_n$.
  
  Consider the element (where we write the columns of tableaux using exponential notation to 
  indicate the multiplicity of each column)
   \begin{align*}
     w_1 := \left(
    \begin{array}{c}
       r \\ \vdots \\2 \\ 1
    \end{array}
     \right)^s \otimes \left(
     \begin{array}{c}
       n+1 \\ \vdots \\ n + 3 - r\\ 1
     \end{array}
     \right)^s \; \in C_1.
   \end{align*}
 Under $\prA$, $w_1$ maps to 
 \begin{align*}
   w_2 := \left(
   \begin{array}{c}
     r+1 \\ r \\ \vdots \\ 3 \\ 2
   \end{array}
   \right)^s \otimes \left(
   \begin{array}{c}
     n+1 \\ \vdots \\ n + 4 - r\\ 2 \\1
   \end{array}
   \right)^s.
  \end{align*}
  However it is not hard to see that
  \begin{align*}
    \Psi(w_1) = \left(
    \begin{array}{c}
      r \\ r-1 \\ \vdots \\ 1
    \end{array}
    \right)^{s-1}
    \otimes
    \left(
    \begin{array}{c}
      r \\ r-1 \\ \vdots \\ 1
    \end{array}
    \right)
    \left(
    \begin{array}{c}
      n+1 \\ \vdots \\ n + 3 - r \\ 1
    \end{array}
    \right)^{s} \in B( (s-1)\omega_r ) \otimes B( (s+1) \omega_r ).
  \end{align*}
  Similarly,
  \begin{align*} 
    \Psi(w_2) = \left(
    \begin{array}{c}
      r+1 \\ r\\ \vdots \\ 5 \\ 4 \\ 2
    \end{array}
    \right)^s \otimes \left(
    \begin{array}{c}
      n+1 \\ \vdots \\ n+4 -r \\ 3\\ 2\\ 1
    \end{array}
    \right)^s \in B(s\omega_{r-1} ) \otimes B( s\omega_{r+1} ) .
  \end{align*}
  Since $\prA'$ preserves the components under $\Psi$ by definition, it cannot agree 
  with $\prA$ on $w_1$. Hence $\prA\neq \prA'$ on $C_1$.

  For the restriction on $C_n$, consider the element
  \begin{align*}
     w_n := \left(
     \begin{array}{c}
       n \\ r - 1 \\ \vdots \\ 1
     \end{array}
     \right)^s \otimes \left(
     \begin{array}{c}
       n+1 \\ n \\ \vdots \\ n + 2 - r
     \end{array}
     \right)^s \in C_n.
   \end{align*}
   Note that $\Psi(w_n)$ is in the $B( (s-1)\omega_r ) \otimes B( (s+1) \omega_r )$
   component.  On the other hand, the image of $w_n$ under $\prA$ is 
   \begin{align*}
     w_{n+1} := \left(
     \begin{array}{c}
        n+1 \\ r \\ \vdots \\ 2
     \end{array}
      \right)^s \otimes \left(
     \begin{array}{c}
        n+1 \\ \vdots \\ n + 3 - r \\ 1
     \end{array}
     \right)^s,
   \end{align*}
   which is in the $B( s\omega_{r-1} ) \otimes B( s\omega_{r+1} )$ component under $\Psi$. So we
   conclude that $\prA \neq \prA'$ on $C_n$.
\end{proof}

For the next proposition, $B:= B(s'\omega_{r'}) \otimes B(s \omega_r)$ is an $A_n$ crystal, 
with $n\ge r' \ge r \ge 1$, and $D:= B(s'\omega_{r'-1}) \otimes B(s\omega_{r-1})$ is an $A_{n-1}$ crystal.
Here we interpret $\omega_0$ as the zero weight.

\begin{proposition} \label{Prop:Induction}
 If $\prA_D$ is the only promotion operator defined on $D$, then any promotion operator $\pr$
 on $B$ agrees with $\prA_B$ on both $C_1(B)$ and $C_n(B)$.
 If $\prA_D$ and $\prA_D'$ are the only promotion operators defined on $D$, then any promotion
 operator $\pr$ on $B$ either agrees with $\prA_B$ on both $C_1(B)$ and $C_n(B)$ or it agrees
 with $\prA_B'$ on both $C_1(B)$ and $C_n(B)$.
\end{proposition}
\begin{proof}
  For the purpose of this proof we set $C_i:=C_i(B)$. Note that while the graphs $C_1$ and
  $C_{n+1}$ are twisted-isomorphic to $D$, the other graphs $C_i$ for $i\neq 1,n+1$ do not 
  have enough arrows.

  We claim there is a unique isomorphism from the $A_{n-1}$ crystal $C_{n+1}$ to $D$.  We define one
  such isomorphism $\phi_{n+1}$, by simply removing the top row from each factor. (It is easy to
  see that this is removing all letters $n+1$).  By the decomposition~\eqref{eq:classicalDecomposition},
  this isomorphism is unique.

  Similary, there is a unique twisted isomorphism from the $\{2, \ldots, n\}$ crystal $C_1$ to
  the $\{1, \ldots, n-1\}$ crystal $D$ with twist given by $\tau: i \mapsto i-1$.
  We define one such isomorphism, $\phi_1$, by removing the bottom row from each factor and
  subtracting one from each letter.  Now
  given any $\tau$-twisted isomorphism $\phi: C_1 \to D$, we get $\phi^{-1} \circ \phi_1: D \to D$
  is a crystal isomorphism.  By~\eqref{eq:classicalDecomposition} this must be the identity, which 
  implies that $\phi = \phi_1$.
  
  Note that this implies that there is a unique twisted isomorphism from $C_{n+1} \to C_1$ with
  twist given by $\tau: i \mapsto i+1$.  This is given by $\phi_{1}^{-1} \circ \phi_{n+1}$, and any
  other $\tau$-twisted isomorphism would give a nontrivial automorphism of $D$.  Since any 
  promotion operator $\pr$ on $B$ must give such a $\tau$-twisted isomorphism when restricted 
  to $C_{n+1}$, we have shown that the promotion operator $\pr$ restricted to $\pr: C_{n+1} \to C_1$ 
  is determined.
  
  Now let $\pr$ be any promotion operator on $B$, restricted to the union of the sets
  $C_1,\ldots,C_{n+1}$.  Let $\prA_B$ be the standard promotion operator on $B$, and $\prA_D$ be the
  standard promotion operator on $D$. Define a map
  $\phi_2 : C_2 \to D$ by
  \begin{align*}
    \phi_2 := \prA_D \circ \phi_1 \circ \prA_B^{-1}.
  \end{align*}
  (All functions written here will be acting on the left.)

  Define a map from $D$ to itself by
  \begin{align*}
    \rho := \phi_2 \circ \pr \circ \phi_1^{-1}.
  \end{align*}
  Note that $\rho$ is a map which takes $D$ to itself and affects content and arrows according to
  axioms~\ref{pt.content} and~\ref{pt.e} of promotion operators as defined in 
  Definition~\ref{definition.promotion}.  We will show that $\rho$ also satisfies the axiom that
  $\rho^n = \id$, hence proving that $\rho$ is a promotion operator on $D$.
  
  We determine the order of $\rho$ by constructing a commutative diagram:
  $$\xymatrix{
  C_1 \ar[d]_{\phi_1} \ar[r]^\pr & C_2 \ar[d]_{\phi_2} \ar[r]^\pr &
    C_3 \ar[d]_{\phi_3} \ar[r]^\pr & {\dots} \ar[r]^\pr 
    & C_{n+1} \ar[r]^{\pr} \ar[d]_{\phi_{n+1}} & C_1 \ar[dl]_{\phi_1} \\
  D \ar[r]^{\rho}              & D \ar[r]^{\rho}              &
    D \ar[r]^{\rho}              & {\dots} \ar[r]^{\rho} 
    & D &
  }
 $$
 The maps $\phi_i$ in this diagram (for $3 \le i \le n$) are defined by requiring this diagram to
 commute.  Specifically, we have 
 \begin{align*}
    \phi_i := \rho \circ \phi_{i-1} \circ \pr^{-1}.
 \end{align*}
 Notice that, by the uniqueness of $\phi_{n+1}$, we must have $\rho \circ \phi_n \circ
 \pr^{-1} = \phi_{n+1}$ on $C_{n+1}$.  Since $\pr$ is a promotion operator, the composition
 $\pr^{n+1}$ along the top row of this diagram must be equal to the identity.  Thus we can collapse
 the diagram to
 $$\xymatrix{
 C_1 \ar[d]_{\phi_1} \ar[dr]^{\phi_1} & \\
 D \ar[r]^{\rho^n} & D
 }
 $$
 which implies $\rho^n = \id$.  This completes the proof that $\rho$ is a promotion operator.

 Now assume that the only choice for a promotion operator on $D$ is $\prA_D$. Recall
 \begin{align*}
   \rho &= \phi_2 \circ \pr \circ \phi_1^{-1}\\
        &= \left( \prA_D \circ \phi_1 \circ \prA_B^{-1} \right) \circ \pr \circ \phi_1^{-1}.\\
   \intertext{Since $\rho = \prA_D$ we multiply both sides by $\prA_D^{-1}$ on the left to obtain}
   \id_D &= \phi_1 \circ \prA_B^{-1} \circ \pr \circ \phi_1^{-1}.\\
   \intertext{Conjugating by $\phi_1$ gives}
   \id_{C_1} &= \prA_B^{-1} \circ \pr\\
   \intertext{which implies}
   \pr &= \prA_B \quad \text{on $C_1$.}
 \end{align*}
 
 Next assume that there are two choices for the promotion operator on $D$, namely $\prA_D$
 and $\prA_D'$. If $\rho=\prA_D$, then by the same arguments as above, we conclude that
 $\pr=\prA_B$ on $C_1$. If $\rho=\prA_D'$, then $\pr=  \phi_2^{-1} \circ \rho \circ \phi_1$ must be
 different from $\prA_B$ on $C_1$. Furthermore, there are no more than these two possibilities
 for $\pr$ on $C_1$. By Lemma~\ref{lemma.prAandprA'}, $\prA_B$ and $\prA_B'$ are different
 on $C_1$. Hence if $\rho=\prA_D'$, then $\pr=\prA_B'$ on $C_1$.

 We now wish to show that $\pr^{-1} = \prA_B^{-1}$ on $C_{n+1}$ if the only choice for a promotion
 operator on $D$ is $\prA_D$.  We keep our definitions of $\phi_1$ and $\phi_{n+1}$ and define 
 $\phi_n: C_n \to D$ by
 \begin{align*}
   \phi_n := \prA_D^{-1} \circ \phi_{n+1} \circ \prA_B.
 \end{align*}
 We now redefine $\rho: D \to D$ (though it will in fact coincide with the old definition) by
 \begin{align*}
   \rho := \phi_{n+1} \circ \pr \circ \phi_n^{-1}.
 \end{align*}
 We again conclude that $\rho$ satisfies the content and arrow properties of a promotion operator,
 and we determine its order with the following diagram:
  $$\xymatrix{
  C_{n+1} \ar[d]_{\phi_{n+1}} & C_n \ar[d]_{\phi_n} \ar[l]_\pr &
  C_{n-1} \ar[d]_{\phi_{n-1}} \ar[l]_\pr & {\dots} \ar[l]_\pr &
  C_1 \ar[l]_{\pr} \ar[d]_{\phi_1} & C_{n+1} \ar[dl]^{\phi_{n+1}} \ar[l]_{\pr} \\
  D               & D \ar[l]^{\rho}              &
    D \ar[l]^{\rho}              & {\dots} \ar[l]^{\rho} 
    & D \ar[l]^{\rho}&
  }
 $$
 (Again the undefined vertical arrows are defined soley to make the diagram commute.) As before, we
 conclude that $\rho^n = \id$ and hence $\pr = \prA_B$ on $C_n$. By very similar arguments as
 before, $\pr$ on $C_n$ is either $\prA_B$ and $\prA_B'$ if there are the two choices $\prA_D$
 and $\prA_D'$ for promotion operators on $D$.
 
 If $\prA_D$ is the only promotion operator on $D$, then we are done. If $\prA_D$ and 
 $\prA_D'$ are the two choices for promotion operators on $D$, then it remains to show
 that $\pr$ is $\prA_B$ on both $C_1$ and $C_n$ or that $\pr$ is $\prA_B'$ on both $C_1$
 and $C_n$.  Recall $w_1$ and $w_n$ as defined in the proof of Lemma~\ref{lemma.prAandprA'}.
 Note that $w_1$ is related to $w_n$ by a series of $f_j$ operators (not including $f_n$).
 Thus $\pr(w_1)$ determines $\pr(w_n)$ and conversely.
\end{proof}

\subsection{Rectangle tensor rectangle case}
\label{subsection.rectangle}

In this section we prove Theorem~\ref{theorem.main} for $(s,r)\neq (s',r')$, $r' \ge r \ge 1$, $r' > 1$, and 
$n\le r+r'-1$.  The proof is by induction on $n$, showing that there is a unique promotion on the $A_n$
crystal $B = B(s'\omega_{r'}) \otimes B(s \omega_r)$, assuming that the statement is true by
induction for smaller $n$.
 
\begin{lemma} \label{Lem:Smalln}
  For $n < r + r' - 2$, the promotion operator on the $A_n$ crystal 
  $B := B(s'\omega_{r'}) \otimes B(s \omega_r)$ is determined.
\end{lemma}
\begin{proof}
   It suffices to show that the promotion operator on $\{1, \ldots, n-1\}$ highest weight elements is
   determined.  The right tensor factor of such an element $w$ must be of the form
   \begin{equation} \label{eq:rightHwFactor}
     \left(
     \begin{array}{c}
       r \\ r-1 \\ \vdots \\ 1
     \end{array}
     \right)^a
     \left(
     \begin{array}{c}
       n+1 \\ r-1 \\ \vdots \\ 1
     \end{array}
     \right)^b.
   \end{equation}
   Hence the bottom row of the left tensor factor can only contain the letters $1$, $r$ and $r+1$
   (since $r'>1$ the letter $n+1$ is not possible in the first row). But if $r$ or $r+1$ appears in the 
   bottom of a column in the left tensor factor, then $r+r'-1$ is the smallest possible number which 
   could appear at the top of that column by columnstrictness.  But since $n<r+r'-2$, 
   this letter is not in our crystal. Hence the bottom row of the left tensor factor consists only
   of $1$s, so $w$ is in $C_1(B)$ and by Proposition~\ref{Prop:Induction} promotion is given by
   either $\prA$ or $\prA'$.
\end{proof}

For the rest of the proof, we assume that $n=r'+r-2$ or $n=r' + r -1$.  The unique element of weight
$s\omega_r$ in $B(s\omega_r)$ is called \emph{Yamanouchi}; it is the tableau with row $i$ filled
with letter $i$ for $1\le i\le r$. After a preliminary lemma, we first show that promotion on all
$\{1,2,\ldots,n-1\}$ highest weight vectors, where the right factor is Yamanouchi, is determined.
Then we prove the claim for general $\{1,2,\ldots,n-1\}$ highest weight elements.

\begin{lemma} \label{Lem:RecRecBracketing}
  Suppose $w$ is a $\{1, \ldots, n-1\}$ highest weight element whose right factor is Yamanouchi
  or whose right factor has height one.  
  Let $\pr, \pr'$ be any two (\emph{a priori} different) promotion operators.
  Let $w_1, w_n$ be defined by
  \begin{equation*}
  \begin{split}
	&w_n \overset{\pr'}{\to} w \overset{\pr'}{\to} w_1.\\
        \text{If} \quad  &w_n \prto w_0 \prto w_1 \quad \text{then} \quad w_0 = w.
  \end{split}
  \end{equation*}
\end{lemma}
\begin{proof}
   We first assume that the right factor has height one.  We claim that $w$ is completely
   specified by:
   \begin{enumerate}
     \item The fact that $w$ is $\{1, \ldots, n-1\}$ highest weight;
     \item The content of $w$; 
     \item The content of the right factor of $w$.
   \end{enumerate}
   Suppose $w$ is $\{1, \ldots, n-1\}$ highest weight.  Then the right factor of $w$ must be of the
   form $1^a (n+1)^b$.  Suppose also that $m_{n+1}(w) = b + c$.  The $A_{n-1}$ crystal consisting of
   those elements of $B$ with $c$ copies of $n+1$ in the left factor and $b$ copies of $n+1$ in the
   right factor is isomorphic to $B( (s'-c) \omega_{r'} + c \omega_{r'-1} ) \otimes B(
   (s-b)\omega_1)$.  This has a decomposition into classical components according to the
   multiplication of Schur functions $s_{(s'^{r-1}, s'-c)} s_{(s-b)}$.  Since this product is
   indexed by a ``near rectangle'' and a rectangle, Theorem~2.1 of~\cite{Stembridge.2001} gives that
   this product, when expanded into Schur functions, is multiplicity free.  Hence, there is at most
   one highest weight vector of a given content.  This proves the claim.

   It is clear that (1) and (2) can be reconstructed from $w_1$.  From (1), we know the right factor
   of $w$ is of the form $1^a (n+1)^b$, and the right factor of $w_0$ is of the form $1^{a'}
   (n+1)^{b'}$.  So it suffices show that $b' = b$.  For this, we note that the $\{2, \dots, n\}$
   lowest weight element associated to $w$ has precisely $s' + b$ copies of $n+1$.  Thus the $\{1,
   \dots, n-1\}$ lowest weight element associated to $w_n$ has $s' + b$ copies of $n$, and the $\{2,
   \dots, n\}$ lowest weight element corresponding to $w_0$ must have $s' + b$ of copies of $n+1$.
   But this is only the case if $b' = b$.

   Now we consider the case that $r > 1$ and the right factor of $w$ is Yamanouchi.  We claim that
   $w$ is completely specified by:
   \begin{enumerate}
     \item The fact that $w$ is $\{1, \ldots, n-1\}$ highest-weight;
     \item The content of $w$; 
     \item The fact that the right factor of $w$ is Yamanouchi.
   \end{enumerate}
   In this case, we know that all $c:= m_{n+1}(w)$ copies of $n+1$ occur in the top row of the left
   factor.  Thus the $A_{n-1}$ crystal consisting of all elements with $c$ copies of $n+1$ in the
   left factor and none in the right is isomorphic to $B( (s'-c) \omega_{r'} + c \omega_{r'-1})
   \otimes B(s \omega_r)$.  Again, the corresponding product of Schur functions is indexed by a near
   rectangle and a rectangle, and so the product is multiplicity free.  Hence there is at most one
   highest weight vector of a given content.

   It is clear that (1) and (2) can be reconstructed from $w_1$.  From (1), we know the right factor
   of $w_0$ must be of the form~\eqref{eq:rightHwFactor}.  We must show that $b = 0$
   in~\eqref{eq:rightHwFactor}.
   
   Assume first that $r' > r$. Let $m_r$ be the number of letters $r$ in $w$.  We note that the number of
   letters $r$ in the $\{2, \ldots, n\}$ highest weight associated to $w$ must also 
   be precisely $m_r$; the only elements that can change at all are the letters $n+1$, and because 
   they are in a row of height $ >r$, they cannot become letters $r$.  Thus the number of letters $r-1$ 
   in the $\{1, \ldots, n-1 \}$ highest
   weight associated to $w_n$ is also $m_r$.  From this we conclude that in any $w_0$, the number of
   $r$s in the associated $\{2, \dots, n\}$ highest weight is $m_r$.  If there were an $n+1$ in the
   right factor of $w$, this would not be true, so we can conclude that $b = 0$.

   Finally suppose that $r' = r$.  Let $k \ge 0$ be the number of columns of the left factor of $w$ 
   whose top two entries are of the form $\genfrac{}{}{0pt}{}{n+1}{r-1}$.   If $k > 0$, the top two rows of 
   the left factor of $w$ is of the general form: 
   \begin{equation*}
   \left(
   \begin{array}{ccccccccc} 
     \le n & \dots & \le n  & n+1 & \dots & n+1 & n+1 & \dots & n+1 \\
     r-1  & \dots & r-1 & r-1 & \dots & r-1 & \ge r &\dots & \ge r\\
     &&&& \dots &&&&
   \end{array}
   \right).
   \end{equation*}
   If $k \le s$, we note that the number of $r$s in the $\{2, \dots, n\}$ highest weight associated to $w$ 
   is again $m_r$.  Hence we can repeat the $r' > r$ argument above to conlude that $b = 0$.  If
   $k > s$, we set $w'$ to be the $\{r, \dots, n\}$ lowest weight associated to $w$.  We have
   $m_{n+1}(w') = s'$, since every element in the top row of the left factor can be raised to an
   $n+1$, and every $r$ in the right factor will be 'blocked' by at least one of the $k$ letters $(n+1)$
   above an $r-1$ on the left on its way up.  Translating this property to $w_n$ and back to $w_0$, we 
   see that the $\{r, \dots, n\}$ lowest weight associated to any $w_0$ must contain the letter $n+1$ 
   precisely $s'$ times.  But the number of $(n+1)$s in this lowest weight must be at least $s' + b$; 
   hence $b= 0$.
   \end{proof}

\begin{lemma} \label{Lem:YamanouchiHighWeights}
   The promotion operator is determined on the set of $\{1, \ldots, n-1\}$ highest weight elements 
   in $B$ for which the right factor is Yamanouchi.
\end{lemma}
\begin{proof}
 For $n=r'+r-2$, the bottom row of the left factor can only contain the letter $1$ by
 the same arguments as in the proof of Lemma~\ref{Lem:Smalln}. Hence all $\{1, \ldots, n-1\}$
 highest weight elements are in $C_1(B)$ and by Proposition~\ref{Prop:Induction} the statement
 follows by induction.
 
 For $n=r'+r-1$, let $w$ be a $\{1, \ldots, n-1\}$ highest weight element with a Yamanouchi
 right factor. The general strategy is to consider the $\prA$-orbit of $w$:
 $$
 w \to w_1 \to \cdots \to w_n \to w.
 $$
 We show that $w_1, \ldots, w_{r'}$ are associated to $\{1,\ldots,n-1\}$ highest weight elements in
 $C_1(B)$ whose promotion we know by induction. We then show that promotion inverse of $w_{r'+1},
 \ldots, w_n$ is determined by showing that the associated $\{2,\ldots,n+1\}$ lowest weight elements
 are in $C_{n+1}(B)$. By Lemma~\ref{Lem:RecRecBracketing} we also know that $w$ is the unique
 element which is simultaneously $\pr(w_n)$ and $\pr^{-1}(w_1)$.

 The right factor of $w$ is Yamanouchi and hence of the form
 \begin{equation*}
   \left(
   \begin{array}{c}
     r  \\ \vdots \\ 1
   \end{array}
   \right)^{s}.
 \end{equation*}
 Note also that the bottom row of the left factor of any element of this crystal must consist of
 letters $\le r+1$; if there is a larger letter in the first row, this would force a letter larger
 than $n+1 = r + r'$ in the top row. In particular, if $r+1$ appears in the bottom row of the left factor,
 then the column has an $n+1$ on the top. Now consider the element $w_i := \prA^i(w)$, for 
 $1 \le i \le r'$.  The right factor of $w_i$ is
 \begin{equation*}
   \left(
   \begin{array}{c}
     r + i \\  \vdots \\ 1 + i
   \end{array}
   \right)^s.
 \end{equation*}
 The bottom row of the left factor contains only letters $< r + i$, which are not bracketed with the 
 right factor.  In particular, we see that the bottom row can always be lowered via $e_j$ operators 
 (without $e_n$) to a row of $1$s.  Since the bottom row of the right factor can also always be lowered 
 to $1$s (without using $e_n$), we see that $w_i$ can be lowered to $C_1(B)$ and so the promotion 
 $w_i \prto w_{i+1}$ for $1 \le i \le r'$ is determined by Proposition~\ref{Prop:Induction}.
 
 Now notice that for $r' \le i \le n$, the top row of the right factor of $w_i$ consists of only
 $(n+1)$s.  Thus the associated $\{2, \ldots, n\}$ lowest weight element is in $C_{n+1}(B)$, and so by
 induction Proposition~\ref{Prop:Induction} we have determined promotion inverse.  Hence we have
 computed $w_1 \prto w_2 \prto \cdots \prto w_n$ and so by Lemma~\ref{Lem:RecRecBracketing} we know
 the orbit of $w$.
 \end{proof}
 
\begin{lemma} \label{Lem:SimpleHighWeights}
   The promotion operator is determined on the set of $\{1, \ldots, n-1\}$ highest weight elements 
   in $B$ with a right factor of
    \begin{equation} \label{eq:hwRightFactor1}
   \left(
   \begin{array}{c}
     n+1 \\ r - 1 \\ r-2 \\ \vdots \\ 1
   \end{array}
   \right)^s.
   \end{equation}
\end{lemma}
\begin{proof}
 As before we set $w_i=\prA^i(w)$.  Every element in $\{w_1, w_2, \ldots, w_{n+1-r}\}$ is associated
 via a sequence of $e_j$ (not including $e_n$) to a $\{1, \ldots, n-1\}$ highest weight element with
 a Yamanouchi right factor. So promotion of these elements is determined by
 Lemma~\ref{Lem:YamanouchiHighWeights}. The remaining elements $\{w_{n+2 -r},\ldots,w_n, w\}$ are
 associated with a $\{2, \ldots, n\}$ lowest weight element in $C_{n+1}(B)$, so promotion inverse is
 determined by Proposition~\ref{Prop:Induction}. Thus the orbit of $w$ is determined. 
\end{proof}

\begin{lemma} \label{Lem:HighWeights}
   The promotion operator is determined on the set of $\{1, \ldots, n-1\}$ highest weight elements 
   in $B$.
\end{lemma}
\begin{proof}
 If $w$ is a $\{1, \ldots, n-1\}$ highest weight element and has a right factor which is Yamanouchi or 
 of the form~\eqref{eq:hwRightFactor1}, the result follows from 
 Lemmas~\ref{Lem:YamanouchiHighWeights} and~\ref{Lem:SimpleHighWeights}.
 Hence we may assume that the top row of the right factor of $w$ contains both the letters $r$ and
 $n+1$. Then the letters in the top row of the right factor of $w_i := \prA^i(w)$ are given by:
 \begin{align*}
   &0: (r,n+1) && 1: (r,r+1) && \dots && n-r: (n-1,n) \\
   & n-r+1: (n, n+1) &&n-r+2: (n+1, n+1) && \dots && n: (n+1, n+1).
 \end{align*}
 Notice that the right factor of $w_i$ for $1\le i<n-r+1$ can be transformed to the Yamanouchi
 element using a sequence of $e_j$ (not including $e_n$). Hence by 
 Lemma~\ref{Lem:YamanouchiHighWeights} promotion on this element is known.  
 In the case that $r=1$, we have determined $w_1 \prto \dots \prto w_n$, and hence by 
 Lemma~\ref{Lem:RecRecBracketing} we have determined the orbit of $w$.  If $r > 1$, then the 
 top row of the right factor of $w_i$ for $n-r+1<i\le n$ consists only of $n+1$, and hence a sequence 
 of $e_j$ (not including $e_n$) can transform these $w_i$ into a $\{1,2,\ldots,n-1\}$ highest 
 weight element with right factor of the form~\eqref{eq:hwRightFactor1}, whose promotion orbit is 
 already determined. In $w_{n-r+1}$, the right factor has the form
 \begin{align*}
   \left(
   \begin{array}{c}
     n \\ \vdots \\ n - r + 1
   \end{array}
   \right)^b
   \left(
   \begin{array}{c}
     n+1 \\ \vdots \\ n-r + 2
   \end{array}
   \right)^a.
 \end{align*}
 Notice that every letter in the right factor is fully bracketed except for the letters $n$.
 Thus every letter $1,2,\ldots,n$ in the left factor of $w_{n-r+1}$ is unbracketed (with respect to
 the right factor).  In particular, every letter in the first row of the left factor is unbracketed,
 so we can reduce them to $1$s. This gives an element of $C_1(B)$ as the associated $\{1, \ldots,
 n-1\}$ highest weight element, and hence promotion is known by induction by
 Proposition~\ref{Prop:Induction}.  Thus the orbit of $w$ is determined, and this completes the
 proof.  
 \end{proof}

By Lemma~\ref{Lem:HighWeights} promotion on all $\{1,2,\ldots,n-1\}$ highest weight elements
is determined. Hence by the Highest Weight Lemma~\ref{lemma.highestWeight} promotion is 
determined on all of $B$. This concludes the induction step in the proof of Theorem~\ref{theorem.main}
when $(s,r)\neq (s',r')$.

\subsection{Equal Tensor Factors} 
\label{subsection.equal_factors}

  Let $B:=B(s\omega_{r}) \otimes B(s\omega_r)$ be the tensor product of two identical classical
  highest weight crystals of type $A_n$ with $n\ge 2$ and $r>1$.  We show in this
  section that there are two promotion operators on $B$, given by the connected operator $\prA$ and
  the disconnected operator $\prA'=\Psi^{-1}\circ\prA\circ\Psi$.
   
  By Proposition~\ref{Prop:Induction}, there are at most two possibilities for the action of promotion on
  the subsets of $B$ given by $C_1:=C_1(B)$ and $C_n:=C_n(B)$.  If promotion restricted to these 
  subsets is given by $\prA$, then all the arguments from Section~\ref{subsection.rectangle}
  apply as before and we are done.  So for the rest of this section we consider the case where 
  promotion on $C_1$ and $C_n$ is given by $\prA'$.  As before by the Highest Weight 
  Lemma~\ref{lemma.highestWeight}, it suffices to determine promotion on all $\{1, \dots, n-1\}$ 
  highest weight elements.

  \begin{lemma} \label{Lem:ef_yam}
   Suppose $\pr$ on $B$ coincides with $\prA'$ on $C_1$ and $C_n$.
   If $w \in B$ is a $\{1,\ldots,n-1\}$ highest weight element, with $\Psi(w) \in B_1 := B((s-1) \omega_r)
   \otimes B((s+1) \omega_r)$, and the right factor of $\Psi(w)$ is Yamanouchi, then the orbit of
   $w$ is given by $\prA'$.
  \end{lemma}
  
  \begin{proof}
  We first note that the conditions of the lemma imply that the right factor of $w$ is Yamanouchi: Suppose 
  $\Psi(w):=v_1 \otimes v_2 \in B_1$ is $\{1, \ldots, n-1\}$ highest weight with $v_2$ being Yamanouchi.
  Then every letter $n+1$ in $v_1.v_2$ is in a row higher than row $r$. Furthermore, since $\Psi$ is a
  crystal isomorphism, $v_1\otimes v_2$ is also $\{1, \dots, n-1\}$ highest weight. Now let 
  $v_1' \otimes v_2' := \Psi^{-1}(v_1 \otimes v_2)$ (so $v_1'.v_2' = v_1.v_2$). This must still be 
  $\{1, \dots, n-1\}$ highest weight, and thus the right tensor factor must be of the 
  form~\eqref{eq:rightHwFactor}. However, any $n+1$ in $v_2'$ would certainly give an $n+1$ at 
  height $r$ in $v_1'.v_2'$.  Thus the only possibility for $v_1'.v_2'$ to agree with $v_1.v_2$ is if $v_2'$ 
  is Yamanouchi.

  Now, we label the elements of the orbit of $w$ under $\prA'$ by
  \begin{displaymath}
    w \to w_1' \to \dots \to w_n' \to w.
  \end{displaymath}
  Recall that $\Psi$ is a crystal isomorphism and hence commutes with the crystal operators and 
  preserves content.  In particular, $\Psi^{-1}(C_1(B_1)) \subset C_1(B)$.  By the proof of 
  Lemma~\ref{Lem:YamanouchiHighWeights} we know that, for $1 \le i \le r$, $\Psi(w_i')$ is connected 
  to $C_1(B_1)$ by a series of classical crystal operators (not involving $e_n$).  Hence $w_i'$ is
   connected to $C_1(B)$.  From the same lemma, it also follows that for $r \le i \le n$, $\Psi(w_i')$ 
   is connected to $C_n(B_1)$ by a series of classical crystal operators (not involving $f_n$); 
   hence $w_i'$ is connected to $C_n(B)$.  Thus the partial orbit
   \begin{displaymath}
     w_1' \to w_2' \to \dots \to w_n'
   \end{displaymath}
   is determined. By Lemma~\ref{Lem:RecRecBracketing}, the entire orbit is now determined.
\end{proof}

  \begin{lemma} \label{Lem:ef_notyam}
     Suppose $\pr$ is a promotion operator on $B$ which coincides with $\prA'$ on $C_1$ and 
     $C_n$. If $w \in B$ is such that $\Psi(w) \in B_1$, then $\pr(w) = \prA'(w)$.
  \end{lemma}
  \begin{proof}
    It remains to show this for those $\{1,\dots,n-1\}$ highest weight elements whose image under 
    $\Psi$ is in $B_1$ and does not have a Yamanouchi right factor.  First consider those elements 
    $w$ where $\Psi(w)$ has only a single repeated column on the right.  Again, we label the orbit under
    $\prA'$ of $w$ by
    \begin{displaymath}
      w_0 := w \to w_1' \to \dots \to w_n' \to w_0.
    \end{displaymath}
    By the proof of Lemma~\ref{Lem:SimpleHighWeights}, $\Psi(w_i')$ for $0 \le i \le n$ is connected 
    by special sequences of crystal operators to elements whose promotion is already determined. 
    Thus this is also true for $w_i'$.  In particular, promotion of $w$ is determined.  This logic can also
    be applied to the remaining $\{1, \ldots, n-1\}$ highest weight elements under
    consideration following the proof of Lemma~\ref{Lem:HighWeights}.
  \end{proof}

  The fact that $\pr$ agrees with $\prA'$ on $B_1$ implies that $\pr(B_2) = B_2$, where
  $B_2 := B(s \omega_{r-1}) \otimes B(s \omega_{r +1 })$.  By Section~\ref{subsection.rectangle}
  we already know that promotion on a tensor product of two distinct rectangles is given by $\prA$; 
  thus we have in this case that $\pr = \prA$ on $B_2$ and thus $\pr = \prA'$ on $B$.

\section{Evidence for Conjecture~\ref{conjecture.promotion}}
\label{section.evidence}
In this section, we present evidence for Conjecture~\ref{conjecture.promotion}.  
In Section~\ref{subsection.factorization} we present theoretical results that support the claims
of the conjecture and in Section~\ref{subsection.computer} we discuss computer evidence.
  
\subsection{Unique factorization into rectangular Schur functions}
\label{subsection.factorization}

We have seen in Lemma~\ref{lemma.existence} that $\prA$ is a valid
promotion operator on a tensor product of classical highest weight
crystals of type $A_n$ indexed by rectangles; furthermore $\prA$
yields a connected affine crystal.

In the trailer of this section, we further argue that two distinct
tensor products of classical highest weight crystals of type $A_n$
indexed by rectangles have non isomorphic classical structures, as
desired for Conjecture~\ref{conjecture.promotion} (otherwise, the two
associated promotion operators could induce two non isomorphic
connected affine crystals).
This statement translates as follow at the level of symmetric
polynomials.
\begin{proposition}
\label{proposition.unique_factorization}
    Let $n\geq 1$. If a symmetric polynomial $P:=P(x_1,\ldots,x_{n+1})$ can be factored as a
    product $P=s_{(c_1^{r_1})}\cdots s_{(c_k^{r_k})}$ of nontrivial rectangular Schur polynomials
    with $1\le r_i\le n$, then this is the unique factorization of $P$ as a product of rectangular
    Schur polynomials.
\end{proposition}

This turns out to be a special case of the following theorem.
\begin{theorem}[Rajan~\cite{Rajan.2004}] \label{theorem.uniqueness}
Let $\mathfrak{g}$ be any simple Lie algebra, and $V_1,\ldots,V_n$ and $W_1,\ldots,W_m$ 
be nontrivial, finite-dimensional, irreducible $\mathfrak{g}$-modules.
If $V_1\otimes \cdots \otimes V_n \cong W_1 \otimes \cdots \otimes W_m$, then $n=m$ and
$V_i \cong W_{\tau(i)}$ for some permutation $\tau$.
\end{theorem}
In type $A$, Purbhoo and van Willigenburg~\cite{Purbhoo.Willigenburg.2007} give a 
combinatorial proof for products of two arbitrary Schur functions. The following combinatorial proof
of Proposition~\ref{proposition.unique_factorization} handles products
of an arbitrary number of rectangular Schur functions.

\begin{proof}[Proof of Proposition~\ref{proposition.unique_factorization}]
    We may impose a total order on rectangular partitions by defining $(c^r) \ge ({c'}^{r'})$ if $r > r'$ or 
    $r = r'$ and $c \ge c'$.  We show that the factor in $P$ indexed by the largest rectangle in this order 
    is uniquely determined. Hence induction on the largest factor proves the proposition.
    
    Without loss of generality we may assume that $(c_1^{r_1}),\ldots,(c_k^{r_k})$ are ordered
    in weakly decreasing order. We use two facts, easily derived from the Littlewood--Richardson 
    rule. Let $Q = \prod_{i = 1}^{k} s_{\lambda^{(i)}}$ be any product of Schur functions.  Let
    $(\nu^{(j)})_{j=1}^m$ be the list of partitions which index the expansion of $Q$ into the
    sum of Schur functions (the order of this list does not matter).  Then
    \begin{enumerate}
      \item For all pairs $(i,j)$ with $1 \le i \le k$ and $1 \le j \le m$ the diagram of
        $\nu^{(j)}$ contains the diagram of $\lambda^{(i)}$.
      \item If $\mu$ is a diagram consisting of the $\lambda^{(i)}$ concatenated to form a partition
        shape, then $\mu$ is one of the $\nu^{(j)}$. 
    \end{enumerate}
    Using these properties, we shall see that $(c_1^{r_1})$, defined to be the index of rectangle
    $\lambda^{(1)}$, can be determined from the collection of $\nu^{(j)}$.   We first find $r_1$.
    Note that property~(1) implies that the height of every diagram $\nu^{(j)}$ is at least $r_1$.
    But by property~(2), there is some shape $\nu^{(j)}$ consisting of the shapes $\lambda^{(i)}$
    concatenated from left to right.  In particular, this shape has height exactly equal to
    $r_1$.  So $r_1$ is the minimum of the heights of the $\nu^{(i)}$.
    
    Since all other rectangles $(c_i^{r_i})$ for $1\le i\le k$ have height $r_i\le r_1$, we may assume
    without loss of generality that $n=r_1$. Each term $s_\nu$ in the Schur expansion of 
    $P$ can be associated with a highest weight crystal element in $B:=B(c_1^{r_1}) \otimes
    \cdots \otimes B(c_k^{r_k})$ of weight $\nu$. Take the collection of all terms $s_\nu$ in
    the Schur expansion of $P$ such that the first $n-1$ parts of $\nu$ agree with the first
    $n-1$ parts of the partition obtained by concatenating all rectangles $(c_1^{r_1}),\ldots,(c_k^{r_k})$.
    This implies in particular that the corresponding highest weight crystal elements in $B$
    are all Yamanouchi in the first $n-1$ rows. If $c_1,\ldots,c_m$ are the widths of the rectangles
    of height $n=r_1$, then the terms $s_\nu$, with $\nu$ given as above, are in one-to-one
    correspondence with the Schur expansion of the following product of complete symmetric
    functions $h_{c_1} \cdots h_{c_m}$ in two variables.
  
    However, note that for $n=1$ we have
    $h_j(x,1)=\frac{1-x^{j+1}}{1-x}$, so its roots are the nontrivial
    $(j+1)$-th roots of unity. Let us consider two factorizations
    $h_{c_1}\cdots h_{c_m} = h_{c'_1}\cdots h_{c'_{m'}}$,
    and show that they must coincide. Consider the largest $h_j$
    occurring either on the left or the right hand side, and consider
    a primitive $(j+1)$-th root of unity $\Xi$. Then, $\Xi$ is a root of
    $h_{c_i}(1,x)$ for some $i$, and by maximality of $j$, $c_i=j$. We
    can therefore factor out $h_j$ from the left hand side -- and
    similarly from the right hand side -- and apply induction.
\end{proof}

\begin{example}
Let us illustrate the proof of Proposition~\ref{proposition.unique_factorization} in terms
of an example. Take $P=s_{22} s_{11}^2 s_{3}$. In the Schur expansion of $P$
there is a term $s_{74}$, which is obtained by concatenating the four rectangles.
All terms are labeled by partitions with at least two parts. This tells us that the height of the
largest rectangle is $r_1=2$.

To determine the width of the largest rectangle, we consider the highest weight crystal elements that
are Yamanouchi in the first $r_1-1=1$ rows:
\begin{equation*}
\begin{aligned}[1]
    \vcenter{\tab{2,2|1,1}} &\otimes \vcenter{\tab{2|1}} \otimes \vcenter{\tab{2|1}}
    \otimes \vcenter{\tab{1,1,1}}, \quad
    \vcenter{\tab{2,2|1,1}} \otimes \vcenter{\tab{3|1}} \otimes \vcenter{\tab{2|1}}
    \otimes \vcenter{\tab{1,1,1}},\\[2mm]
     \vcenter{\tab{2,3|1,1}} &\otimes \vcenter{\tab{2|1}} \otimes \vcenter{\tab{2|1}}
    \otimes \vcenter{\tab{1,1,1}}, \quad
    \vcenter{\tab{3,3|1,1}} \otimes \vcenter{\tab{2|1}} \otimes \vcenter{\tab{2|1}}
    \otimes \vcenter{\tab{1,1,1}}.
\end{aligned}
\end{equation*}
The second row of these elements gives precisely the expansion of the $\mathfrak{sl}_2$
or two variable expansion of $h_2 h_1^2$, which is unique. Hence $c_1=2$, and the largest rectangle 
is $(2^2)$.
\end{example}

\subsection{Computer exploration}
\label{subsection.computer}
The research was partially driven by computer exploration. In
particular, we implemented a branch-and-bound algorithm to search for
all (connected) (weak) promotion operators on a given classical
crystal. The algorithm goes down a search tree, deciding progressively
to which $\{2,\dots,n\}$-component each connected
$\{1,\dots,n-1\}$-component is mapped by promotion. Branches are cut
as soon as it can be decided that the yet partially defined promotion
cannot satisfy condition~(2') of Remark~\ref{remark.shift2}, or cannot
be connected. The algorithm can also take advantage of the symmetries
of the classical crystal (not fully though, by lack of appropriate
group theoretical tools in \texttt{MuPAD}), and uses some heuristics
for the decision order.  The branch cutting works reasonably well; for
the difficult case of $B(1) ^{\otimes4}$ in type $A_2$, where the
total search space is a priori of size $144473849856000$, with
$2!3!3!=72$ symmetries, the algorithm actually explores $115193$
branches in $5$ hours and $26$ minutes (on a $2$ GHz Linux PC), using
$16$M of memory. The result is $8$ isomorphic connected promotion
operators: $9$ symmetries out of the $72$ could be exploited to cut
the search space.

\begin{example}
  We start by loading the \texttt{MuPAD-Combinat} package, and setting
  the notation for tensor products:.
\begin{Mexin}
package("MuPAD-Combinat"):  
operators::setTensorSymbol("#"):
\end{Mexin}

Consider the $A_2$ classical crystal $C:=B(1)\otimes B(1)\otimes B(1)$:
\begin{Mexin}
B1 := crystals::tableaux(["A",2], Shape = [1]):
C := B1 # B1 # B1:
\end{Mexin}
The decomposition into classical components is given by $s_1^3 = s_3 +
s_{111} + 2 s_{21}$ (note the multiplicity of $s_{21}$). There are four
promotion operators:
\begin{Mexin}
promotions := C::promotions():
nops(promotions)
\end{Mexin}
\begin{Mexout}
                                      4
\end{Mexout}
Let us construct the corresponding crystal graphs:
\begin{Mexin}
  affineCrystals :=
  [crystals::affineFromClassicalAndPromotion(C, promotion)
   $ promotion in promotions]:
\end{Mexin}
Among them, two are connected:
\begin{Mexin}
  [ A::isConnected() $ A in affineCrystals ]
\end{Mexin}
\begin{Mexout}
                          [TRUE, FALSE, FALSE, TRUE]
\end{Mexout}
But they are in fact isomorphic via the exchange of the two $(2,1)$-classical components:
\begin{Mexin}
nops((affineCrystals[1])::isomorphisms(affineCrystals[4]))
\end{Mexin}
\begin{Mexout}
                                      1
\end{Mexout}
The other two affine crystals are disconnected, and induced by the
decomposition $s_1^3 = s_1s_{11} + s_1s_2$ of $s_1^3$ into a sum of
products of rectangles. The use of the options \texttt{Connected} and
\texttt{Symmetries} cuts down the search tree. It turns out that for
our current crystal, the symmetries are fully exploited, and we only
get one isomorphic copy of the connected promotion operator:
\begin{Mexin}
nops(C::promotions(Connected, Symmetries))
\end{Mexin}
\begin{Mexout}
                                      1
\end{Mexout}

Now consider the $A_2$ classical crystal $C:=B(2,1)\otimes B(2,1)$.
\begin{Mexin}
B21 := crystals::tableaux(["A",2], Shape = [2,1]):
C := B21 # B21:
\end{Mexin}
The highest weights of the classical crystal are given by the
following Schur polynomial expansion:
\begin{equation}
    s_{21}^2 = s_{42} + s_{411} + s_{33} + 2s_{321} + s_{222}.
\end{equation}
Beware that, since $n=2$, the term $s_{2211}$ is zero. Also, the
crystal for $s_{411}$ is isomorphic to that for $s_3$, and
similarly $s_{222}$ is trivial. Finally, note the multiplicity of $s_{321}$.

There are no connected promotion operators:
\begin{Mexin}
nops(C::promotions(Connected))
\end{Mexin}
\begin{Mexout}
                                       0
\end{Mexout}
Indeed, $f$ has no factorization into products of rectangle
Schur polynomials. On the other hand, there are eight disconnected
promotion operators:
\begin{Mexin}
nops(C::promotions())
\end{Mexin}
\begin{Mexout}
                                       8
\end{Mexout}
They are induced by the following four decompositions:
\begin{equation}
    \begin{aligned}
      s_{21}^2 
      &= s_{22} s_{1} s_{1} + s_{3}s_{111}\\
      &= s_{22} s_{11} + s_{22} s_{2} + s_{3}s_{111}\\
      &= s_{22} s_{2} + s_{2} s_{1}s_{111} + s_{33}\\
      &= s_{11} s_{11} s_{2} + s_{33}
    \end{aligned}
\end{equation}
combined with the automorphism which exchanges the two
$(3,2,1)$-classical components.
\end{example}

The examples of Figure~\ref{fig:promotionsForA11} for $n=1$ were found
with this algorithm. On the other hand, we ran systematic tests on the
following crystals with $n\geq 2$:
\begin{itemize}
\item All tensor products of rows with up to $3$ factors and up to $6$
  cells (except $B(2)^{\otimes 3}$) in type $A_2$ and up to $7$ cells
  (except $B(3)\otimes B(2)^{\otimes 2}$) in type $A_3$ and $A_4$;
\item $B(3,2,1)\otimes B(1)$, $B(2,1)\otimes B(2,1)$, $B(2,1)\otimes
  B(1)\otimes B(1)$, $B(2,2)\otimes B(1,1,1)$, $B(2,2)\otimes
  B(1,1,1)\otimes B(1)$, $B(1)^{\otimes 4}$, in type $A_2$ and $A_3$.
\end{itemize}
They all agree with Conjectures~\ref{conjecture.promotion}
and~\ref{conjecture.Kashiwara}. Namely, for tensor products of
rectangles, there is a unique connected promotion operator, up to
isomorphism; for other tensor products, there is none.

In the smaller examples, we further checked that the total number of
promotions was exactly the number automorphism of the underlying
classical crystal (that is $\prod m_{\lambda}!$ where $m_\lambda$ is
the number of classical components of highest weight $\lambda$) times
the number of decompositions of the symmetric function $\sum m_\lambda
s_\lambda$ into sums of products of rectangular Schur functions.

\bibliographystyle{alpha}
\bibliography{promotion}
\end{document}

%% file: figure-promotions-B11B13-1.tex
\begin{tikzpicture}[>=latex,join=bevel,scale=.5]
\tiny%
  \node (N_1) at (50bp,306bp) [draw,draw=none] {$1 \!\otimes\!111$};
  \node (N_2) at (5bp,232bp) [draw,draw=none] {$1 \!\otimes\!112$};
  \node (N_3) at (5bp,158bp) [draw,draw=none] {$1 \!\otimes\!122$};
  \node (N_4) at (95bp,232bp) [draw,draw=none] {$2 \!\otimes\!111$};
  \node (N_5) at (95bp,158bp) [draw,draw=none] {$2 \!\otimes\!112$};
  \node (N_6) at (5bp,84bp) [draw,draw=none] {$1 \!\otimes\!222$};
  \node (N_7) at (95bp,84bp) [draw,draw=none] {$2 \!\otimes\!122$};
\node (N_8) at (50bp,10bp) [draw,draw=none] {$2 \!\otimes\!222$};
  \draw [->,DarkBlue] (N_6) to [bend right=10] node [left] {$1$} (N_8);
  \draw [<-] (N_6) to [bend left=10] node [right] {$0$} (N_8);
  \draw [->,DarkBlue] (N_1) to [bend right=10] node [left] {$1$} (N_2);
  \draw [<-] (N_1) to [bend left=10] node [right] {$0$} (N_2);
  \draw [->,DarkBlue] (N_3) to [bend right=10] node [left] {$1$} (N_6);
  \draw [<-] (N_3) to [bend left=10] node [right] {$0$} (N_6);
  \draw [->,DarkBlue] (N_4) to [bend right=10] node [left] {$1$} (N_5);
  \draw [<-] (N_4) to [bend left=10] node [right] {$0$} (N_5);
  \draw [->,DarkBlue] (N_5) to [bend right=10] node [left] {$1$} (N_7);
  \draw [<-] (N_5) to [bend left=10] node [right] {$0$} (N_7);
  \draw [->,DarkBlue] (N_2) to [bend right=10] node [left] {$1$} (N_3);
  \draw [<-] (N_2) to [bend left=10] node [right] {$0$} (N_3);
\end{tikzpicture}

%% file: figure-promotions-B11B13-2.tex
\begin{tikzpicture}[>=latex,join=bevel,scale=.5]
\tiny%
  \node (N_1) at (50bp,306bp) [draw,draw=none] {$1 \!\otimes\!111$};
  \node (N_2) at (5bp,232bp) [draw,draw=none] {$1 \!\otimes\!112$};
  \node (N_3) at (5bp,158bp) [draw,draw=none] {$1 \!\otimes\!122$};
  \node (N_4) at (95bp,232bp) [draw,draw=none] {$2 \!\otimes\!111$};
  \node (N_5) at (95bp,158bp) [draw,draw=none] {$2 \!\otimes\!112$};
  \node (N_6) at (5bp,84bp) [draw,draw=none] {$1 \!\otimes\!222$};
  \node (N_7) at (95bp,84bp) [draw,draw=none] {$2 \!\otimes\!122$};
\node (N_8) at (50bp,10bp) [draw,draw=none] {$2 \!\otimes\!222$};

  \draw [<-] (N_5) to [] node [right] {$0$} (N_6);

  \draw [->,DarkBlue] (N_6) to [bend right=10] node [left] {$1$} (N_8);
  \draw [<-] (N_6) to [bend left=10] node [right] {$0$} (N_8);

  \draw [->,DarkBlue] (N_1) to [bend right=10] node [left] {$1$} (N_2);
  \draw [<-] (N_1) to [bend left=10] node [right] {$0$} (N_2);

  \draw [->,DarkBlue] (N_3) to [] node [left] {$1$} (N_6);
  \draw [->,DarkBlue] (N_4) to [] node [left] {$1$} (N_5);
  \draw [<-] (N_3) to [] node [right] {$0$} (N_7);
  \draw [->,DarkBlue] (N_5) to [] node [left] {$1$} (N_7);
  \draw [->,DarkBlue] (N_2) to [] node [left] {$1$} (N_3);
  \draw [<-] (N_2) to [] node [right] {$0$} (N_5);
  \draw [<-] (N_4) to [] node [right] {$0$} (N_3);
\end{tikzpicture}

%% file: figure-promotions-B11B13-3.tex
\begin{tikzpicture}[>=latex,join=bevel,scale=.5]
\tiny%
  \node (N_1) at (50bp,306bp) [draw,draw=none] {$1 \!\otimes\!111$};
  \node (N_2) at (5bp,232bp) [draw,draw=none] {$1 \!\otimes\!112$};
  \node (N_3) at (5bp,158bp) [draw,draw=none] {$1 \!\otimes\!122$};
  \node (N_4) at (95bp,232bp) [draw,draw=none] {$2 \!\otimes\!111$};
  \node (N_5) at (95bp,158bp) [draw,draw=none] {$2 \!\otimes\!112$};
  \node (N_6) at (5bp,84bp) [draw,draw=none] {$1 \!\otimes\!222$};
  \node (N_7) at (95bp,84bp) [draw,draw=none] {$2 \!\otimes\!122$};
\node (N_8) at (50bp,10bp) [draw,draw=none] {$2 \!\otimes\!222$};
  \draw [<-] (N_5) to [] node [right] {$0$} (N_6);
  \draw [<-] (N_7) to [] node [right] {$0$} (N_8);
  \draw [->,DarkBlue] (N_6) to [] node [left] {$1$} (N_8);
  \draw [->,DarkBlue] (N_1) to [] node [left] {$1$} (N_2);
  \draw [->,DarkBlue] (N_3) to [] node [left] {$1$} (N_6);
  \draw [->,DarkBlue] (N_4) to [] node [left] {$1$} (N_5);
  \draw [<-] (N_3) to [] node [right] {$0$} (N_7);
  \draw [->,DarkBlue] (N_5) to [] node [left] {$1$} (N_7);
  \draw [->,DarkBlue] (N_2) to [] node [left] {$1$} (N_3);
  \draw [<-] (N_1) to [] node [right] {$0$} (N_4);
  \draw [<-] (N_2) to [] node [right] {$0$} (N_5);
  \draw [<-] (N_4) to [] node [right] {$0$} (N_3);
\end{tikzpicture}

%% file: figure-promotions-B11B13-4.tex
\begin{tikzpicture}[>=latex,join=bevel,scale=.5]
\tiny%
  \node (N_1) at (50bp,306bp) [draw,draw=none] {$1 \!\otimes\!111$};

  \node (N_2) at (5bp,232bp) [draw,draw=none] {$1 \!\otimes\!112$};
  \node (N_4) at (95bp,232bp) [draw,draw=none] {$2 \!\otimes\!111$};

  \node (N_3) at (5bp,158bp) [draw,draw=none] {$1 \!\otimes\!122$};
  \node (N_5) at (95bp,158bp) [draw,draw=none] {$2 \!\otimes\!112$};

  \node (N_6) at (5bp,84bp) [draw,draw=none] {$1 \!\otimes\!222$};
  \node (N_7) at (95bp,84bp) [draw,draw=none] {$2 \!\otimes\!122$};

\node (N_8) at (50bp,10bp) [draw,draw=none] {$2 \!\otimes\!222$};
  \draw [<-] (N_7) to [] node [right] {$0$} (N_8);
  \draw [->,DarkBlue] (N_6) to [] node [left] {$1$} (N_8);
  \draw [->,DarkBlue] (N_1) to [] node [left] {$1$} (N_2);

  \draw [->,DarkBlue] (N_3) to [bend right=10] node [left] {$1$} (N_6);
  \draw [<-] (N_3) to [bend left=10] node [right] {$0$} (N_6);

  \draw [->,DarkBlue] (N_4) to [bend right=10] node [left] {$1$} (N_5);
  \draw [<-] (N_4) to [bend left=10] node [right] {$0$} (N_5);

  \draw [->,DarkBlue] (N_5) to [bend right=10] node [left] {$1$} (N_7);
  \draw [<-] (N_5) to [bend left=10] node [right] {$0$} (N_7);

  \draw [->,DarkBlue] (N_2) to [bend right=10] node [left] {$1$} (N_3);
  \draw [<-] (N_2) to [bend left=10] node [right] {$0$} (N_3);

  \draw [<-] (N_1) to [] node [right] {$0$} (N_4);
\end{tikzpicture}

%% file: figure-row3xrow2-inversionLess.tex
\begin{tikzpicture}
  \tikzstyle{content}=[draw,rectangle,fill=gray!25]
  \tikzstyle{bracket}=[draw,rectangle]
  \tikzstyle{twopath}=[draw=gray!25,rectangle]
  \node (11111) at (0bp,250bp) [content] {$11\!\otimes\!111$};
  \node (11112) at (0bp,200bp) [twopath] {$11\!\otimes\!112$};
  \node (11122) at (0bp,150bp) [twopath] {$11\!\otimes\!122$};
  \node (11113) at (70bp,200bp) [twopath] {$11\!\otimes\!113$};
  \node (11123) at (70bp,150bp) [twopath] {$11\!\otimes\!123$};
  \node (11223) at (70bp,100bp) [twopath] {$11\!\otimes\!223$};
  \node (11133) at (140bp,150bp) [twopath] {$11\!\otimes\!133$};
  \node (12223) at (70bp,50bp) [bracket] {$12\!\otimes\!223$};
  \node (11233) at (140bp,100bp) [bracket] {$11\!\otimes\!233$};
  \node (11222) at (0bp,100bp) [twopath] {$11\!\otimes\!222$};
  \node (22223) at (70bp,0bp) [twopath] {$22\!\otimes\!223$};
  \node (12233) at (140bp,50bp) [twopath] {$12\!\otimes\!233$};
  \node (11333) at (210bp,100bp) [twopath] {$11\!\otimes\!333$};
  \node (12222) at (0bp,50bp) [twopath] {$12\!\otimes\!222$};
  \node (22233) at (140bp,0bp) [twopath] {$22\!\otimes\!233$};
  \node (12333) at (210bp,50bp) [twopath] {$12\!\otimes\!333$};
  \node (22222) at (0bp,0bp) [content] {$22\!\otimes\!222$};
  \node (22333) at (210bp,0bp) [twopath] {$22\!\otimes\!333$};
  \node (13333) at (280bp,50bp) [twopath] {$13\!\otimes\!333$};
  \node (23333) at (280bp,0bp) [twopath] {$23\!\otimes\!333$};
  \node (33333) at (350bp,0bp) [content] {$33\!\otimes\!333$};
\draw [->,DarkBlue] (11111) to node [auto] {$1$} (11112);
\draw [->,DarkBlue] (11112) to node [auto] {$1$} (11122);
\draw [->,DarkRed] (11112) to node [auto] {$2$} (11113);
\draw [->,DarkBlue] (11122) to node [auto] {$1$} (11222);
\draw [->,DarkRed] (11122) to node [auto] {$2$} (11123);
\draw [->,DarkBlue] (11113) to node [auto] {$1$} (11123);
\draw [->,DarkBlue] (11123) to node [auto] {$1$} (11223);
\draw [->,DarkRed] (11123) to node [auto] {$2$} (11133);
\draw [->,DarkBlue] (11223) to node [auto] {$1$} (12223);
\draw [->,DarkRed] (11223) to node [auto] {$2$} (11233);
\draw [->,DarkBlue] (11133) to node [auto] {$1$} (11233);
\draw [->,DarkBlue] (12223) to node [auto] {$1$} (22223);
\draw [->,DarkRed] (12223) to node [auto] {$2$} (12233);
\draw [->,DarkBlue] (11233) to node [auto] {$1$} (12233);
\draw [->,DarkRed] (11233) to node [auto] {$2$} (11333);
\draw [->,DarkBlue] (11222) to node [auto] {$1$} (12222);
\draw [->,DarkRed] (11222) to node [auto] {$2$} (11223);
\draw [->,DarkRed] (22223) to node [auto] {$2$} (22233);
\draw [->,DarkBlue] (12233) to node [auto] {$1$} (22233);
\draw [->,DarkRed] (12233) to node [auto] {$2$} (12333);
\draw [->,DarkBlue] (11333) to node [auto] {$1$} (12333);
\draw [->,DarkBlue] (12222) to node [auto] {$1$} (22222);
\draw [->,DarkRed] (12222) to node [auto] {$2$} (12223);
\draw [->,DarkRed] (22233) to node [auto] {$2$} (22333);
\draw [->,DarkBlue] (12333) to node [auto] {$1$} (22333);
\draw [->,DarkRed] (12333) to node [auto] {$2$} (13333);
\draw [->,DarkRed] (22222) to node [auto] {$2$} (22223);
\draw [->,DarkRed] (22333) to node [auto] {$2$} (23333);
\draw [->,DarkBlue] (13333) to node [auto] {$1$} (23333);
\draw [->,DarkRed] (23333) to node [auto] {$2$} (33333);
\end{tikzpicture}

%% file: figure-row3xrow3-inversionLess.tex
\begin{tikzpicture}
  \tikzstyle{content}=[draw,rectangle,fill=gray!25]
  \tikzstyle{bracket}=[draw,rectangle]
  \tikzstyle{twopath}=[draw=gray!25,rectangle]
  \node (111111) at (0bp,300bp) [content] {$111\!\otimes\!111$};
  \node (111112) at (0bp,250bp) [twopath] {$111\!\otimes\!112$};
  \node (111122) at (0bp,200bp)  [twopath] {$111\!\otimes\!122$};
  \node (111113) at (80bp,250bp)  [twopath] {$111\!\otimes\!113$};
  \node (111123) at (80bp,200bp)  [twopath] {$111\!\otimes\!123$};
  \node (111223) at (80bp,150bp)  {$111\!\otimes\!223$};
  \node (111133) at (160bp,200bp)  [twopath] {$111\!\otimes\!133$};
  \node (112223) at (80bp,100bp)  [twopath] {$112\!\otimes\!223$};
  \node (111233) at (160bp,150bp)  {$111\!\otimes\!233$};
  \node (111222) at (0bp,150bp)  {$111\!\otimes\!222$};
  \node (122223) at (80bp,50bp) [bracket] {$122\!\otimes\!223$};
  \node (112233) at (160bp,100bp) [bracket] {$112\!\otimes\!233$};
  \node (111333) at (240bp,150bp)  {$111\!\otimes\!333$};
  \node (112222) at (0bp,100bp)  [twopath] {$112\!\otimes\!222$};
  \node (222223) at (80bp,0bp)  [twopath] {$222\!\otimes\!223$};
  \node (122233) at (160bp,50bp)  [twopath] {$122\!\otimes\!233$};
  \node (112333) at (240bp,100bp)  {$112\!\otimes\!333$};
  \node (122222) at (0bp,50bp) [twopath] {$122\!\otimes\!222$};
  \node (222233) at (160bp,0bp)  [twopath] {$222\!\otimes\!233$};
  \node (122333) at (240bp,50bp)  {$122\!\otimes\!333$};
  \node (113333) at (320bp,100bp)  [twopath] {$113\!\otimes\!333$};
  \node (222222) at (0bp,0bp) [content] {$222\!\otimes\!222$};
  \node (222333) at (240bp,0bp)  {$222\!\otimes\!333$};
  \node (123333) at (320bp,50bp)  [twopath] {$123\!\otimes\!333$};
  \node (223333) at (320bp,0bp)  [twopath] {$223\!\otimes\!333$};
  \node (233333) at (400bp,0bp)  [twopath] {$233\!\otimes\!333$};
  \node (333333) at (480bp,0bp) [content] {$333\!\otimes\!333$};
  \node (133333) at (400bp,50bp) [twopath] {$133\!\otimes\!333$};
\draw [->,DarkBlue] (111111) to node [auto] {$1$} (111112);
\draw [->,DarkBlue] (111112) to node [auto] {$1$} (111122);
\draw [->,DarkRed] (111112) to node [auto] {$2$} (111113);
\draw [->,DarkBlue] (111122) to node [auto] {$1$} (111222);
\draw [->,DarkRed] (111122) to node [auto] {$2$} (111123);
\draw [->,DarkBlue] (111113) to node [auto] {$1$} (111123);
\draw [->,DarkBlue] (111123) to node [auto] {$1$} (111223);
\draw [->,DarkRed] (111123) to node [auto] {$2$} (111133);
\draw [->,DarkBlue] (111223) to node [auto] {$1$} (112223);
\draw [->,DarkRed] (111223) to node [auto] {$2$} (111233);
\draw [->,DarkBlue] (111133) to node [auto] {$1$} (111233);
\draw [->,DarkBlue] (112223) to node [auto] {$1$} (122223);
\draw [->,DarkRed] (112223) to node [auto] {$2$} (112233);
\draw [->,DarkBlue] (111233) to node [auto] {$1$} (112233);
\draw [->,DarkRed] (111233) to node [auto] {$2$} (111333);
\draw [->,DarkBlue] (111222) to node [auto] {$1$} (112222);
\draw [->,DarkRed] (111222) to node [auto] {$2$} (111223);
\draw [->,DarkBlue] (122223) to node [auto] {$1$} (222223);
\draw [->,DarkRed] (122223) to node [auto] {$2$} (122233);
\draw [->,DarkBlue] (112233) to node [auto] {$1$} (122233);
\draw [->,DarkRed] (112233) to node [auto] {$2$} (112333);
\draw [->,DarkBlue] (111333) to node [auto] {$1$} (112333);
\draw [->,DarkBlue] (112222) to node [auto] {$1$} (122222);
\draw [->,DarkRed] (112222) to node [auto] {$2$} (112223);
\draw [->,DarkRed] (222223) to node [auto] {$2$} (222233);
\draw [->,DarkBlue] (122233) to node [auto] {$1$} (222233);
\draw [->,DarkRed] (122233) to node [auto] {$2$} (122333);
\draw [->,DarkBlue] (112333) to node [auto] {$1$} (122333);
\draw [->,DarkRed] (112333) to node [auto] {$2$} (113333);
\draw [->,DarkBlue] (122222) to node [auto] {$1$} (222222);
\draw [->,DarkRed] (122222) to node [auto] {$2$} (122223);
\draw [->,DarkRed] (222233) to node [auto] {$2$} (222333);
\draw [->,DarkBlue] (122333) to node [auto] {$1$} (222333);
\draw [->,DarkRed] (122333) to node [auto] {$2$} (123333);
\draw [->,DarkBlue] (113333) to node [auto] {$1$} (123333);
\draw [->,DarkRed] (222222) to node [auto] {$2$} (222223);
\draw [->,DarkRed] (222333) to node [auto] {$2$} (223333);
\draw [->,DarkBlue] (123333) to node [auto] {$1$} (223333);
\draw [->,DarkRed] (123333) to node [auto] {$2$} (133333);
\draw [->,DarkRed] (223333) to node [auto] {$2$} (233333);
\draw [->,DarkRed] (233333) to node [auto] {$2$} (333333);
\draw [->,DarkBlue] (133333) to node [auto] {$1$} (233333);
\end{tikzpicture}


%% file: figure-row4xrow2-inversionLess.tex
\begin{tikzpicture}
  \tikzstyle{content}=[draw,rectangle,fill=gray!25]
  \tikzstyle{bracket}=[draw,rectangle]
  \tikzstyle{twopath}=[draw=gray!25,rectangle]
  \node (111111) at (0bp,300bp) [content] {$11\!\otimes\!1111$};
  \node (111112) at (0bp,250bp) [twopath] {$11\!\otimes\!1112$};
  \node (111122) at (0bp,200bp) [twopath] {$11\!\otimes\!1122$};
  \node (111113) at (80bp,250bp) [twopath] {$11\!\otimes\!1113$};
  \node (111123) at (80bp,200bp) [twopath] {$11\!\otimes\!1123$};
  \node (111223) at (80bp,150bp) [twopath] {$11\!\otimes\!1223$};
  \node (111133) at (160bp,200bp) [twopath] {$11\!\otimes\!1133$};
  \node (112223) at (80bp,100bp) [twopath] {$11\!\otimes\!2223$};
  \node (111233) at (160bp,150bp) [twopath] {$11\!\otimes\!1233$};
  \node (111222) at (0bp,150bp) [twopath] {$11\!\otimes\!1222$};
  \node (122223) at (80bp,50bp) [bracket] {$12\!\otimes\!2223$};
  \node (112233) at (160bp,100bp) [bracket] {$11\!\otimes\!2233$};
  \node (111333) at (240bp,150bp) [twopath] {$11\!\otimes\!1333$};
  \node (112222) at (0bp,100bp) [twopath] {$11\!\otimes\!2222$};
  \node (222223) at (80bp,0bp) [twopath] {$22\!\otimes\!2223$};
  \node (122233) at (160bp,50bp) [twopath] {$12\!\otimes\!2233$};
  \node (112333) at (240bp,100bp) [bracket] {$11\!\otimes\!2333$};
  \node (122222) at (0bp,50bp) [twopath] {$12\!\otimes\!2222$};
  \node (222233) at (160bp,0bp) [twopath] {$22\!\otimes\!2233$};
  \node (122333) at (240bp,50bp) [twopath] {$12\!\otimes\!2333$};
  \node (113333) at (320bp,100bp) [twopath] {$11\!\otimes\!3333$};
  \node (222222) at (0bp,0bp) [content] {$22\!\otimes\!2222$};
  \node (222333) at (240bp,0bp) [twopath] {$22\!\otimes\!2333$};
  \node (123333) at (320bp,50bp) [twopath] {$12\!\otimes\!3333$};
  \node (223333) at (320bp,0bp) [twopath] {$22\!\otimes\!3333$};
  \node (233333) at (400bp,0bp) [twopath] {$23\!\otimes\!3333$};
  \node (333333) at (480bp,0bp) [content] {$33\!\otimes\!3333$};
  \node (133333) at (400bp,50bp) [twopath] {$13\!\otimes\!3333$};
\draw [->,DarkBlue] (111111) to node [auto] {$1$} (111112);
\draw [->,DarkBlue] (111112) to node [auto] {$1$} (111122);
\draw [->,DarkRed] (111112) to node [auto] {$2$} (111113);
\draw [->,DarkBlue] (111122) to node [auto] {$1$} (111222);
\draw [->,DarkRed] (111122) to node [auto] {$2$} (111123);
\draw [->,DarkBlue] (111113) to node [auto] {$1$} (111123);
\draw [->,DarkBlue] (111123) to node [auto] {$1$} (111223);
\draw [->,DarkRed] (111123) to node [auto] {$2$} (111133);
\draw [->,DarkBlue] (111223) to node [auto] {$1$} (112223);
\draw [->,DarkRed] (111223) to node [auto] {$2$} (111233);
\draw [->,DarkBlue] (111133) to node [auto] {$1$} (111233);
\draw [->,DarkBlue] (112223) to node [auto] {$1$} (122223);
\draw [->,DarkRed] (112223) to node [auto] {$2$} (112233);
\draw [->,DarkBlue] (111233) to node [auto] {$1$} (112233);
\draw [->,DarkRed] (111233) to node [auto] {$2$} (111333);
\draw [->,DarkBlue] (111222) to node [auto] {$1$} (112222);
\draw [->,DarkRed] (111222) to node [auto] {$2$} (111223);
\draw [->,DarkBlue] (122223) to node [auto] {$1$} (222223);
\draw [->,DarkRed] (122223) to node [auto] {$2$} (122233);
\draw [->,DarkBlue] (112233) to node [auto] {$1$} (122233);
\draw [->,DarkRed] (112233) to node [auto] {$2$} (112333);
\draw [->,DarkBlue] (111333) to node [auto] {$1$} (112333);
\draw [->,DarkBlue] (112222) to node [auto] {$1$} (122222);
\draw [->,DarkRed] (112222) to node [auto] {$2$} (112223);
\draw [->,DarkRed] (222223) to node [auto] {$2$} (222233);
\draw [->,DarkBlue] (122233) to node [auto] {$1$} (222233);
\draw [->,DarkRed] (122233) to node [auto] {$2$} (122333);
\draw [->,DarkBlue] (112333) to node [auto] {$1$} (122333);
\draw [->,DarkRed] (112333) to node [auto] {$2$} (113333);
\draw [->,DarkBlue] (122222) to node [auto] {$1$} (222222);
\draw [->,DarkRed] (122222) to node [auto] {$2$} (122223);
\draw [->,DarkRed] (222233) to node [auto] {$2$} (222333);
\draw [->,DarkBlue] (122333) to node [auto] {$1$} (222333);
\draw [->,DarkRed] (122333) to node [auto] {$2$} (123333);
\draw [->,DarkBlue] (113333) to node [auto] {$1$} (123333);
\draw [->,DarkRed] (222222) to node [auto] {$2$} (222223);
\draw [->,DarkRed] (222333) to node [auto] {$2$} (223333);
\draw [->,DarkBlue] (123333) to node [auto] {$1$} (223333);
\draw [->,DarkRed] (123333) to node [auto] {$2$} (133333);
\draw [->,DarkRed] (223333) to node [auto] {$2$} (233333);
\draw [->,DarkRed] (233333) to node [auto] {$2$} (333333);
\draw [->,DarkBlue] (133333) to node [auto] {$1$} (233333);
\end{tikzpicture}


%% file: figure-row5xrow1-inversionLess.tex
\begin{tikzpicture}
  \tikzstyle{content}=[draw,rectangle,fill=gray!25]
  \tikzstyle{bracket}=[draw,rectangle]
  \tikzstyle{twopath}=[draw=gray!25,rectangle]
  \node (111111) at (0bp,300bp) [content] {$1\!\otimes\!11111$};
  \node (111112) at (0bp,250bp) [twopath] {$1\!\otimes\!11112$};
  \node (111122) at (0bp,200bp) [twopath] {$1\!\otimes\!11122$};
  \node (111113) at (80bp,250bp) [twopath] {$1\!\otimes\!11113$};
  \node (111123) at (80bp,200bp) [twopath] {$1\!\otimes\!11123$};
  \node (111223) at (80bp,150bp) [twopath] {$1\!\otimes\!11223$};
  \node (111133) at (160bp,200bp) [twopath] {$1\!\otimes\!11133$};
  \node (112223) at (80bp,100bp) [twopath] {$1\!\otimes\!12223$};
  \node (111233) at (160bp,150bp) [twopath] {$1\!\otimes\!11233$};
  \node (111222) at (0bp,150bp) [twopath] {$1\!\otimes\!11222$};
  \node (122223) at (80bp,50bp) [bracket] {$1\!\otimes\!22223$};
  \node (112233) at (160bp,100bp) [twopath] {$1\!\otimes\!12233$};
  \node (111333) at (240bp,150bp) [twopath] {$1\!\otimes\!11333$};
  \node (112222) at (0bp,100bp) [twopath] {$1\!\otimes\!12222$};
  \node (222223) at (80bp,0bp) [twopath] {$2\!\otimes\!22223$};
  \node (122233) at (160bp,50bp) [bracket] {$1\!\otimes\!22233$};
  \node (112333) at (240bp,100bp) [twopath] {$1\!\otimes\!12333$};
  \node (122222) at (0bp,50bp) [twopath] {$1\!\otimes\!22222$};
  \node (222233) at (160bp,0bp) [twopath] {$2\!\otimes\!22233$};
  \node (122333) at (240bp,50bp) [bracket] {$1\!\otimes\!22333$};
  \node (113333) at (320bp,100bp) [twopath] {$1\!\otimes\!13333$};
  \node (222222) at (0bp,0bp) [content] {$2\!\otimes\!22222$};
  \node (222333) at (240bp,0bp) [twopath] {$2\!\otimes\!22333$};
  \node (123333) at (320bp,50bp) [bracket] {$1\!\otimes\!23333$};
  \node (223333) at (320bp,0bp) [twopath] {$2\!\otimes\!23333$};
  \node (233333) at (400bp,0bp) [twopath] {$2\!\otimes\!33333$};
  \node (333333) at (480bp,0bp) [content] {$3\!\otimes\!33333$};
  \node (133333) at (400bp,50bp) [twopath] {$1\!\otimes\!33333$};
\draw [->,DarkBlue] (111111) to node [auto] {$1$} (111112);
\draw [->,DarkBlue] (111112) to node [auto] {$1$} (111122);
\draw [->,DarkRed] (111112) to node [auto] {$2$} (111113);
\draw [->,DarkBlue] (111122) to node [auto] {$1$} (111222);
\draw [->,DarkRed] (111122) to node [auto] {$2$} (111123);
\draw [->,DarkBlue] (111113) to node [auto] {$1$} (111123);
\draw [->,DarkBlue] (111123) to node [auto] {$1$} (111223);
\draw [->,DarkRed] (111123) to node [auto] {$2$} (111133);
\draw [->,DarkBlue] (111223) to node [auto] {$1$} (112223);
\draw [->,DarkRed] (111223) to node [auto] {$2$} (111233);
\draw [->,DarkBlue] (111133) to node [auto] {$1$} (111233);
\draw [->,DarkBlue] (112223) to node [auto] {$1$} (122223);
\draw [->,DarkRed] (112223) to node [auto] {$2$} (112233);
\draw [->,DarkBlue] (111233) to node [auto] {$1$} (112233);
\draw [->,DarkRed] (111233) to node [auto] {$2$} (111333);
\draw [->,DarkBlue] (111222) to node [auto] {$1$} (112222);
\draw [->,DarkRed] (111222) to node [auto] {$2$} (111223);
\draw [->,DarkBlue] (122223) to node [auto] {$1$} (222223);
\draw [->,DarkRed] (122223) to node [auto] {$2$} (122233);
\draw [->,DarkBlue] (112233) to node [auto] {$1$} (122233);
\draw [->,DarkRed] (112233) to node [auto] {$2$} (112333);
\draw [->,DarkBlue] (111333) to node [auto] {$1$} (112333);
\draw [->,DarkBlue] (112222) to node [auto] {$1$} (122222);
\draw [->,DarkRed] (112222) to node [auto] {$2$} (112223);
\draw [->,DarkRed] (222223) to node [auto] {$2$} (222233);
\draw [->,DarkBlue] (122233) to node [auto] {$1$} (222233);
\draw [->,DarkRed] (122233) to node [auto] {$2$} (122333);
\draw [->,DarkBlue] (112333) to node [auto] {$1$} (122333);
\draw [->,DarkRed] (112333) to node [auto] {$2$} (113333);
\draw [->,DarkBlue] (122222) to node [auto] {$1$} (222222);
\draw [->,DarkRed] (122222) to node [auto] {$2$} (122223);
\draw [->,DarkRed] (222233) to node [auto] {$2$} (222333);
\draw [->,DarkBlue] (122333) to node [auto] {$1$} (222333);
\draw [->,DarkRed] (122333) to node [auto] {$2$} (123333);
\draw [->,DarkBlue] (113333) to node [auto] {$1$} (123333);
\draw [->,DarkRed] (222222) to node [auto] {$2$} (222223);
\draw [->,DarkRed] (222333) to node [auto] {$2$} (223333);
\draw [->,DarkBlue] (123333) to node [auto] {$1$} (223333);
\draw [->,DarkRed] (123333) to node [auto] {$2$} (133333);
\draw [->,DarkRed] (223333) to node [auto] {$2$} (233333);
\draw [->,DarkRed] (233333) to node [auto] {$2$} (333333);
\draw [->,DarkBlue] (133333) to node [auto] {$1$} (233333);
\end{tikzpicture}